\documentclass[10pt]{article}
\usepackage{amsmath}
\numberwithin{equation}{section}
\usepackage{amsfonts}
\usepackage{amssymb}
\usepackage{graphicx}
\usepackage{mathrsfs}
\usepackage{xcolor}
\usepackage{verbatim}
\usepackage{mathrsfs}
\usepackage[body={15.5cm,21cm}, top=3cm]{geometry}
\usepackage{paralist}
\usepackage{ntheorem}
\usepackage{appendix}

\allowdisplaybreaks[4]
\usepackage{hyperref}
\hypersetup{colorlinks=true,
linkcolor=blue,
anchorcolor=blue,
citecolor=blue}
\providecommand{\U}[1]{\protect \rule{.1in}{.1in}}
\newtheorem{theorem}{Theorem}[section]

\newtheorem{definition}[theorem]{Definition}

\newtheorem{lemma}[theorem]{Lemma}

\newtheorem{proposition}[theorem]{Proposition}
\newtheorem{remark}[theorem]{Remark}

\theoremstyle{empty}

\newenvironment{proof}[1][Proof]{\noindent \textbf{#1.} }{\  \rule{0.5em}{0.5em}}
\allowdisplaybreaks[2]
\begin{document}
\title{Stochastic recursive optimal control problem with monotonicity conditions under $G$-framework}
\author{
Qiangjun Tang \thanks{Zhongtai Securities Institute for Financial 
Studies, Shandong University, Jinan, China (tangqj@mail.sdu.edu.cn).}}
\date{}
\maketitle

\begin{abstract}
    In this paper, we study the stochastic
    recursive optimal control problem under non-Lipschitz settings.
    More precisely, we suppose that the generator of $G$-BSDE describing the running cost is 
    uniformly continuous and monotonic with respect to
    the first unknown variable.
    Using the comparison theorem for 
    $G$-BSDE and the stability of viscosity solutions, 
    we establish the dynamic programming principle and
    the connection between
    the value function and the viscosity solution of the 
    associated Hamilton-Jacobi-Bellman equation. We provide an example of continuous time Epstein-Zin utility 
    to demonstrate the application of our study.
\end{abstract}


\textbf{MSC-classification}: 60G65, 60H10, 93E20, 49L25

\section{Introduction}\label{sec1}

In 1992, Duffie and Epstein \cite{DE} pioneered to 
use the solution of BSDE to describe the stochastic differential
utility. Subsequently, the control problem related to 
stochastic differential utility, namely the stochastic 
recursive optimal control problem, was firstly introduced by Peng
\cite{Peng1992}, the dynamic programming principle and the connection
between value function and Hamilton-Jacobi-Bellman (HJB for short) equation
in the Lipschitz setting of the aggregator (or generator from BSDE point of view)
were established. 
Since then, stochastic 
recursive optimal control problem was developed in various settings,
see \cite{Wuyu08, Lipeng09}. Here, we focus on the non-Lipschitz settings, 
that is monotone assumption which was first introduced  
to the classical BSDE theory by Peng \cite{Peng1991} for the infinite horizon BSDE.
Following Peng's footsteps, there are many works adopted 
the monotonic condition to weaken 
the Lipschitz assumption such that BSDE more general, refer to 
\cite{HYPSG95,Pardoux99,PP99,BDHPS03,PZ18,WWZ22}.
Especially, 
Zhuo et al. \cite{ZDP20} consider the stochastic recursive control problem with
the generator of the BSDE is monotonic with respect to
the first unknown variable and uniformly continuous in the second unknown variable.


Recently, Hu et al. \cite{HJPS1,HJb} investigated the following BSDE driven by $G$-Brownian motion:
\begin{small}
    \begin{equation*}
        Y_{t} =\xi
        +\int_{t}^{T}f(s,Y_{s},Z_{s})ds+\sum_{i,j=1}^{d}\int_{t}^{T}g^{ij}(s,Y_{s},Z_{s})d\langle
        B^i,B^j\rangle _{s} 
        -\int_{t}^{T}Z_{s}dB_{s}-(K_{T}-K_{t}).
    \end{equation*}
\end{small}They proved the well-posedness, comparison theorem and Feynman-Kac formula
in the standard Lipschitz settings. Song \cite{Song2019} established existence and uniqueness
with monotonic generators by monotonic approximations (see \cite{HQW20,Liu20,WZ21,He22,Zhao24} for other settings). 
We also mention so-called second-order BSDEs proposed
independently by Soner et al. \cite{STZ12}, which is close to $G$-BSDE.
Possama\"{i} \cite{POD13} extended the results of Soner et al. \cite{STZ12} 
to the case of a generator satisfying a monotonicity condition
in $y$.

Based on Hu et al. \cite{HJPS1,HJb}, 
Hu and Ji \cite{HJ17} studied a stochastic recursive optimal control problem, 
the state equation is governed by the
following controlled SDE driven by a $G$-Brownian motion ($G$-SDE for short)
\begin{small}
    \begin{equation*}
        \left\{
        \begin{array}{l}
            dX_{s}^{t,x,u}=b(s,X_{s}^{t,x,u},u_{s})ds+h^{ij}(s,X_{s}^{t,x
            ,u},u_{s})d\langle B^{i},B^{j}\rangle_{s}+\sigma(s,X_{s}^{t,x,u}
            ,u_{s})dB_{s},\\
            X_{t}^{t,x,u}=x.\\
        \end{array}\right.
    \end{equation*}
\end{small}The cost functional is formulated by the solution to 
the following $G$-BSDE
on the interval $[t,T]$:
\begin{small}
    \begin{align}\label{eq-intro-1}
        Y_{s}^{t,x,u}  &  =\Phi(X_{T}^{t,x,u})+\int_{s}^{T}f(r,X_{r}^{t,x,u},Y_{r}^{t,x,u},Z_{r}^{t,x
        ,u},u_{r})dr
        +\int_{s}^{T}g^{ij}(r,X_{r}^{t,x,u},Y_{r}^{t,x,u},Z_{r}^{t,x,u}
        ,u_{r})d\langle B^{i},B^{j}\rangle_{r}\notag\\
        &\quad  -\int_{s}^{T}Z_{r}^{t,x,u}dB_{r}-(K_{T}^{t,x,u}-K_{s}^{t,x,u}).
    \end{align}
\end{small}And the value function is defined as below
\begin{small}
    \begin{equation*}
        V(t,x):=\underset{u\in\mathcal{U}[t,T]}{ess\inf}\ Y_{t}^{t,x,u}. 
    \end{equation*}
\end{small}They established the dynamic programming principle and 
showed that the value function is
the unique viscosity solution of the related HJB equation under Lipschitz settings.
Other related works that were considered with Lipschitz condition can 
refer to \cite{HW2018,HJL22}. 

We observe that there is still an interesting and challenging problem about the 
dynamic programming principle and related HJB equation 
for stochastic
recursive optimal control problem driven by 
$G$-Brownian motion
with the monotonic aggregator.
Thus, in this paper, we relax the Lipschitz assumptions,
that is the generator of $G$-BSDE \eqref{eq-intro-1} is 
non-Lipschitz, i.e., uniformly continuous, 
monotonic with respect to $y$ and 
uniformly continuous with respect to $u$.

Compare to the standard Lipschitz settings, monotone settings would 
bring many difficulties in the 
dynamic programming principle and related HJB equation 
for stochastic recursive optimal control problem driven by 
$G$-Brownian motion.
Firstly, note that the existence and uniqueness of solutions to $G$-BSDEs under monotonicity 
conditions require the generator to be uniformly continuous with respect to $\omega$ (see \cite{Song2019}). 
However, under assumptions {\bf (B1)}-{\bf (B4)}, the well-posedness of the solution to the state 
equation \eqref{state-2} in the control problem cannot be guaranteed.
Secondly, non-Lipschitz settings bring much trouble
both in the establishment of the dynamic 
programming principle and in the
verification of conditions for the stability of 
viscosity solutions which gives the connection between the value
function and the viscosity solution of the associated HJB equation.
Fortunately, with the help of characterization to spaces $L_G^p(\Omega_T)$ 
and $M_G^p(0,T)$, we only need to slightly modify the results in \cite{Song2019} 
to establish the existence and uniqueness of solutions for the state equation \eqref{state-2}. 
On the other hand, based on the work of Halm 
Lepeltier and San Martin \cite{LSM97},
we adopt a sequence of Lipschitz generators to approximate
our non-Lipschitz generators,
from which we can establish the uniformly convergence of the 
sequences of value function and Hamiltonian in every compact subset of their own
domains.
Therefore, we can establish the dynamic programming principle 
for our stochastic recursive optimal control problem.
Moreover, we show the connection between the value
function and the viscosity solution of the associated HJB equation 
by stability of the viscosity solutions.

In practice,
our stochastic recursive optimal control problem
could be applied to deal with stochastic recursive utility 
with non-Lipschitz setting under $G$-framework.
In $G$-framework, the investor concern the ambiguity or model uncertainty.
Corresponding to our study, we consider the continuous-time Epstein-Zin utility with 
volatility uncertainty, its aggregator has a form 
\begin{small}
    \begin{equation}\label{Epstein-zin-utility}
        f(c,u)=\frac{\delta}{1-\frac{1}{\psi}}(1-\gamma)u
        \Biggl[\biggl(\frac{c_t}{((1-\gamma)u)^{\frac{1}{1-\gamma}}}\biggr)^{1-\frac{1}{\psi}}-1\Biggr],  
\end{equation}
\end{small}where $\delta> 0$ is the rate of time preference, 
$0 <\gamma\neq 1$ is the coefficient of relative risk 
aversion and $0 <\psi\neq 1$
is the elasticity of intertemporal substitution.
We adjust the parameters such $f(c,u)$ in \eqref{Epstein-zin-utility}
is not Lipschitz with
respect $u$ but could be monotonic with respect to $u$. 
Then, we define the value function of the investors by the solution of $G$-BSDE with 
the generator given by \eqref{Epstein-zin-utility}
and indicate the rationale behind this definition. 
Finally, we conclude that the value function of the investor 
is a unique viscosity solution of the HJB equation.

The paper is organized as follows. In Section \ref{sec2}, we recall 
some fundamental results
on $G$-expectation theory. 
Then, we formulate our stochastic recursive optimal control problem 
and establish the dynamic
programming principle in Section \ref{sec4}. In Section \ref{sec5}, 
we derive the fully nonlinear HJB equation
and prove that the value function is 
the unique viscosity solution of the obtained HJB equation.
Finally, we give an example to indicate
the application of our work to finance industry in Section \ref{sec6}.
Some useful estimates are presented in Appendix.

\section{Preliminaries}\label{sec2}

In this section, we review some basic notions and results of $G$-expectation theory and
$G$-BSDEs. More details can be found in \cite{P07a,P10,Song11,DHP11,HJPS1,HJb}. 

\subsection{$G$-expectation}
Let $\Omega_T=C_{0}([0,T];\mathbb{R})$, the space of
real-valued continuous functions with $\omega_0=0$, be endowed
with the supremum norm and 
let  $B$ be the canonical
process. For each $T>0$, set
\begin{small}
    \begin{equation*}
        L_{ip} (\Omega_T):=\{ \varphi(B_{t_{1}},...,B_{t_{n}}):  \ n\in\mathbb {N}, \ t_{1}
        ,\cdots, t_{n}\in[0,T], \ \varphi\in C_{b,Lip}(\mathbb{R}^{ n})\},
    \end{equation*}
\end{small}where $C_{b,Lip}(\mathbb{R}^{ n})$ denotes the set of 
bounded Lipschitz functions on $\mathbb{R}^{n}$.

Denote $\mathbb{S}_d$ be the space of all $d \times d$ symmetric matrices. 
For each given sublinear, continuous 
and monotone function  $G:\mathbb{S}_d\rightarrow\mathbb{R}$,
Peng \cite{P10} constructed a dynamic coherent nonlinear $G$-expectation space
$(\Omega_T, L_{ip} (\Omega_T), \widehat{\mathbb{E}},(\widehat{\mathbb{E}}_t)_{t\ge0})$, 
the corresponding
canonical process $B_t$ is called $G$-Brownian motion. The
function $G:\mathbb{S}_d\to\mathbb{R}$ defined by 
\begin{small}
    \begin{equation}
        G(A):=\frac{1}{2}\widehat{\mathbb{E}}\big[\langle AB_1,B_1\rangle\big], \label{eq-G}
    \end{equation}
\end{small}In this paper we only consider non-degenerate $G$-normal distribution, i.e.,
there exists some $\underline{\sigma}^2 > 0$ such that 
$G(A)-G(B) \ge \underline{\sigma}^2 \text{tr}[A-B]$ 
for any $A \ge B$.

Define $\Vert\xi\Vert_{L_{G}^{p}}:=(\widehat{\mathbb{E}}[|\xi|^{p}])^{1/p}$ 
for $\xi\in L_{ip}(\Omega_T)$ and $p\geq1$.   
The completion of $L_{ip} (\Omega_T)$ under this norm  
is denote by $L_{G}^{p}(\Omega)$. For all $t\in[0,T]$, $\widehat{\mathbb{E}}_t[\cdot]$ 
is a continuous mapping on $L_{ip}(\Omega_T)$ w.r.t the norm $\|\cdot\|_{L_G^1}$. 
Hence, the conditional $G$-expectation $\mathbb{\widehat{E}}_{t}[\cdot]$ can be
extended continuously to the completion $L_{G}^{1}(\Omega_T)$. 
\cite{DHP11} prove that the $G$-expectation has the following representation.
\begin{theorem}[\cite{DHP11}]\label{thm2.1}
	There exists a weakly compact set
	$\mathcal{P}$ of probability
	measures on $(\Omega_T,\mathcal{B}(\Omega_T))$, such that
    \begin{small}
        \begin{equation*}
            \widehat{\mathbb{E}}[\xi]=\sup_{P\in\mathcal{P}}E_{P}[\xi] 
            \text{ for all } \xi\in  {L}_{G}^{1}{(\Omega_T)}.
        \end{equation*}
    \end{small}$\mathcal{P}$ is called a set that represents $\widehat{\mathbb{E}}$.
\end{theorem}

Let $\mathcal{P}$ be a weakly compact set that represents 
$\widehat{\mathbb{E}}$, we define a capacity 
\begin{small}
    \begin{equation*}
        c(A):=\sup_{P\in\mathcal{P}}P(A),\quad A\in\mathcal{B}(\Omega)
    \end{equation*}
\end{small}

A set $A\in\mathcal{B}(\Omega_T)$ is called polar if $c(A)=0$.  A
property holds ``quasi-surely'' (q.s.) if it holds outside a
polar set. In this paper, we do not distinguish two random 
variables $X$ and $Y$ if $X=Y$, q.s.. 

Now we state the monotone convergence theorem of $G$-framework, which is different from
the classical case.
\begin{theorem}[\cite{DHP11}]\label{monotone-con-thm}
    Let $\{X_n\}_{n=1}^\infty\in L_G^1(\Omega_T)$ be such $X_n \downarrow X$ q.s.,
    then $\widehat{\mathbb{E}}[X_n]\downarrow \widehat{\mathbb{E}}[X]$.
    In particular, if $X\in L_G^1(\Omega_T)$, then $\widehat{\mathbb{E}}[Xn-X]\downarrow 0$, as $n\to\infty$.   
\end{theorem} 

\begin{definition}
    A mapping $X$ on $\Omega$ with values in a topological space is said to
    be quasi-continuous (q.c.): if $\forall \epsilon>0$, there exists an 
    open set $O$ with $c(O) < \epsilon$ such that $X|_{O^c}$ is continuous.
\end{definition}

\begin{definition}
    We say that $X : \Omega\mapsto\mathbb{R}$ 
    has a quasi-continuous version 
    if there exists a quasi-continuous function 
    $Y: \Omega\mapsto\mathbb{R}$ such that $X = Y$, q.s..
\end{definition}

\begin{definition}\label{def2.2}
	Let $M_{G}^{0}(0,T)$ be the collection of processes in the
	following form: for a given partition 
    $\{t_{0},\cdot\cdot\cdot,t_{N}\}=\pi_{T}$ of $[0,T]$,
    \begin{small}
        \begin{equation*}
            \eta_{t}(\omega)=\sum_{j=0}^{N-1}\xi_{j}(\omega)\mathbf{1}_{[t_{j},t_{j+1})}(t),
        \end{equation*}
    \end{small}where $\xi_{i}\in L_{ip}(\Omega_{t_{i}})$, $i=0,1,2,\cdot\cdot\cdot,N-1$. 
    For each $p\geq1$ and $\eta\in M_G^0(0,T)$, let 
    $\small \|\eta\|_{H_G^p}:=\tiny \biggl\{\widehat{\mathbb{E}}[(\int_0^T|\eta_s|^2ds)^{p/2}]\biggr\}^{1/p}$, 
    $\small \Vert\eta\Vert_{M_{G}^{p}}:=\tiny \biggl(\mathbb{\widehat{E}}[\int_{0}^{T}|\eta_{s}|^{p}ds]\biggr)^{1/p}$
    and denote by $H_G^p(0,T)$,  $M_{G}^{p}(0,T)$ the completion
	of $M_{G}^{0}(0,T)$ under the norm 
    $\small \|\cdot\|_{H_G^p}$, $\small \|\cdot\|_{M_G^p}$ respectively.
\end{definition}

Let $S_G^0(0,T)=\{h(t,B_{t_1\wedge t}, \ldots,B_{t_n\wedge t}):t_1,\ldots,t_n\in[0,T],h\in C_{b,Lip}(\mathbb{R}^{n+1})\}$. 
For $p\geq 1$ and $\eta\in S_G^0(0,T)$, set 
$\small \|\eta\|_{S_G^p}=\tiny \biggl\{\widehat{\mathbb{E}}\biggl[\sup_{t\in[0,T]}|\eta_t|^p\biggr]\biggr\}^{1/p}$. 
Denote by $S_G^p(0,T)$ the completion of $S_G^0(0,T)$ under the norm $\small \|\cdot\|_{S_G^p}$.  

For each $t\in[0,T]$, denote $L^0(\Omega_t)$
the space of all $\mathcal{B}(\Omega_t)$-measurable real functions.
Then $L_G^p(\Omega_T)$ has the following representation.
\begin{theorem}[\cite{DHP11}]\label{thm-LGp}
    We have
    \begin{small}
        \begin{align*}
            L_G^p(\Omega_T)&=\Big\{X\in L^0(\Omega_T):
            \lim_{N\to\infty}\widehat{\mathbb{E}}\Big[|X|^p\mathbf{1}_{\{|X|\ge N\}}\Big]=0,\text{ and }\\
            &\qquad \quad X \text{ has a quasi-continuous version}\Big\}.
        \end{align*}
    \end{small}   
\end{theorem}

Define, for each $p\ge1$,
\[
    \mathbb{M}^p(0,T)=\Biggl\{\eta: \text{ progressively measurable on } 
    [0,T]\times \Omega_T \text{ and } \widehat{\mathbb{E}}\biggl[
        \int_{0}^{T}|\eta_t|^pdt
    \biggr]<\infty\Biggr\}.
\]
Then $M_G^p(0,T)$ has the following characterization.
\begin{theorem}[\cite{HWZ16}]\label{thm-MGp}
    For each $p\ge1$,
    \begin{align*}
        M_G^p(0,T)&=\Biggl\{
            \eta\in \mathbb{M}^p(0,T):
            \lim_{N\to\infty}\widehat{\mathbb{E}}\Biggl[
                \int_{0}^{T}|\eta_t|^p\mathbf{1}_{\{|\eta_t|\ge N\}}dt
            \Biggr]=0 \text{ and }\\
            &\qquad\quad  \eta \text{ has a quasi-continuous version}
        \Biggr\}.
    \end{align*}
\end{theorem}

The following proposition can be regarded as the 
Burkholder-Davis-Gundy (BDG) inequality under $G$-expectation framework.
\begin{proposition}[\cite{P19}]\label{prop2.3}
	If $\eta\in H_G^{\alpha}(0,T)$ with $\alpha\geq 1$ and $p\in(0,\alpha]$, then we have
    \begin{small}
        \begin{equation*}
            \underline{\sigma}^p c_p\widehat{\mathbb{E}}\biggl[\biggl(\int_0^T |\eta_s|^2ds\biggr)^{p/2}\biggr]\leq
            \widehat{\mathbb{E}}\biggl[\sup_{t\in[0,T]}\biggl|\int_0^t\eta_s dB_s\biggr|^p\biggr]\leq
            \bar{\sigma}^p C_p\widehat{\mathbb{E}}\biggl[\biggl(\int_0^T |\eta_s|^2ds\biggr)^{p/2}\biggr],
        \end{equation*}
    \end{small}where $0<c_p<C_p<\infty$ are constants depending on $p, T$.
\end{proposition}

\begin{definition}\label{def2.4}
A process $\{M_t\}$ with values in $L_G^1(\Omega_T)$ is called 
a $G$-martingale (supermartingale,submartingale) if
$\widehat{\mathbb{E}}_t[M_s]=M_s(\leq,\geq)$ for any $s\leq t$.
\end{definition}

\begin{lemma}[\cite{HJPS1}]\label{lemma2.6}
    Let $X \in S_G^\alpha(0,T)$ for some $\alpha> 1$ and $\alpha^*=\frac{1}{\alpha-1}.$
    Assume that $K^l, l = 1, 2,$ are two non-increasing $G$-martingales with $K^l_0 = 0$ and 
    $K^l_T\in L_{G}^{\alpha^*}(\Omega_T)$. Then the process defined by
    \begin{small}
        \begin{equation*}
            \int_{0}^{t}X_s^+dK_s^1+\int_{0}^{t}X_s^-dK_s^2
        \end{equation*}
    \end{small}is also a non-increasing $G$-martingale.
\end{lemma}

\subsection{$G$-BSDEs with Lipschitz coefficients}
For simplicity, we denote by $\mathfrak{S}_{G}^{\alpha }(0,T)$ the
collection of processes $(Y,Z,K)$ such that $Y\in S_{G}^{\alpha }(0,T)$, 
$Z\in H_{G}^{\alpha }(0,T)$, $K$ is a decreasing $G$-martingale with $K_{0}=0$
and $K_{T}\in L_{G}^{\alpha }(\Omega _{T})$.
Consider the following $G$-BSDEs on $[0, T]$,
\begin{small}
    \begin{align}\label{GBSDE}
        Y_{t} =\xi
        +\int_{t}^{T}f(s,Y_{s},Z_{s})ds+\sum_{i,j=1}^{d}\int_{t}^{T}g^{ij}(s,Y_{s},Z_{s})d\langle
        B^i,B^j\rangle _{s} 
        -\int_{t}^{T}Z_{s}dB_{s}-(K_{T}-K_{t}).
    \end{align}
\end{small} where $f(t,w,y,z), g^{ij}(t,w,y,z):[0,T]\times\Omega\times\mathbb{R}\times\mathbb{R}^d\to\mathbb{R}$ 
satisfy the following:

\begin{enumerate}
    \item [{\bf (H1)}] There exists some $\beta> 1$ such that 
    for any $(y,z)\in\mathbb{R}\times\mathbb{R}^d$, 
    $f(\cdot,\cdot,y,z), g^{ij}(\cdot,\cdot,y,z)\in M_{G}^{\beta}(0,T)$.
    \item [{\bf (H2)}] There exists a constant $L> 0$ such that
    \begin{small}
        \begin{equation*}
            |f(t,y_1,z_1)-f(t,y_2,z_2)|+\sum_{i,j=1}^{d}
            |g^{ij}(t,y_1,z_1)-g^{ij}(t,y_2,z_2)|\le L(|y_1-y_2|+|z_1-z_2|).
        \end{equation*}
    \end{small}
\end{enumerate}

\begin{theorem}[\cite{HJPS1}]\label{thm2.7}
    Suppose that $\xi\in L_{G}^{\beta}(\Omega_T)$ and $f, g^{ij}$ satisfy 
    {\bf (H1)}-{\bf (H2)} for some $\beta>1$. Then, equation \eqref{GBSDE}
    has a unique solution $(Y, Z, K) \in \mathfrak{S}_{G}^{\alpha}(0,T)$
    for any $1<\alpha<\beta$. 
\end{theorem}

\section{Stochastic recursive optimal control problem}\label{sec4}
In this section, we are committed to deriving a generalized dynamic programming 
for the stochastic recursive control problem driven by $G$-Brownian motion.
\subsection{State equations}
We consider admissible controls, which are defined as below.
\begin{definition}
    For each $t\in[0,T]$, $u$ is said to be an admissible control on
    $[t,T]$, if it satisfies the following conditions:
    \begin{enumerate}
    \item [(i)]$u:[t,T]\times\Omega\rightarrow U$ where $U$ is a given compact set of $\mathbb{R}^{m}$;
    \item [(ii)]$u\in M_{G}^{2}(t,T;\mathbb{R}^{m})$.
    \end{enumerate}
\end{definition}
Denote $\mathcal{U}[t,T]$ the set of all admissible controls.
In the following, we use Einstein summation convention.

For each given $t\in[0,T]$, 
$\zeta\in\cup_{\varepsilon>0}L_{G}^{2+\varepsilon}(\Omega_{t};\mathbb{R}^{n})$ 
and $u\in\mathcal{U}[t,T]$, we consider the following $G$-SDE and $G$-BSDE:
\begin{small}
    \begin{equation}\label{state-1}
        \left\{
        \begin{array}{l}
            dX_{s}^{t,\zeta,u}=b(s,X_{s}^{t,\zeta,u},u_{s})ds+h^{ij}(s,X_{s}^{t,\zeta
            ,u},u_{s})d\langle B^{i},B^{j}\rangle_{s}+\sigma(s,X_{s}^{t,\zeta,u}
            ,u_{s})dB_{s},\\
            X_{t}^{t,\zeta,u}=\zeta,\\
        \end{array}\right.
    \end{equation}   
\end{small}and 
\begin{small}
    \begin{align}
        Y_{s}^{t,\zeta,u}  &  =\Phi(X_{T}^{t,\zeta,u})+\int_{s}^{T}f(r,X_{r}^{t,\zeta,u},Y_{r}^{t,\zeta,u},Z_{r}^{t,\zeta
        ,u},u_{r})dr
        +\int_{s}^{T}g^{ij}(r,X_{r}^{t,\zeta,u},Y_{r}^{t,\zeta,u},Z_{r}^{t,\zeta,u}
        ,u_{r})d\langle B^{i},B^{j}\rangle_{r}\notag\\
        &\quad  -\int_{s}^{T}Z_{r}^{t,\zeta,u}dB_{r}-(K_{T}^{t,\zeta,u}-K_{s}^{t,\zeta,u}),\ \ s\in[
            t,T].\label{state-2}
    \end{align}
\end{small}We make the following assumptions on the coefficients of $G$-SDE
\eqref{state-1}.
\begin{enumerate}
    \item [{\bf (A1)}] $b(t,x,u),h^{ij}(t,x,u):[0,T]\times\mathbb{R}^{n}\times U\rightarrow\mathbb{R}^{n}$, 
    $\sigma(t,x,u):[0,T]\times\mathbb{R}^{n}\times U\rightarrow\mathbb{R}^{n\times d}$ 
    are deterministic functions and continuous in $t$;
    \item [{\bf (A2)}] There exists a constant $C>0$ 
    such that for any $t\in [0, T]$, $x,x_1,x_2\in\mathbb{R}^{n}$, $u,v\in U$ 
    \begin{small}
        \begin{align*}
            &  |b(s,x_1,u)-b(s,x_2,v)|+|h^{ij}(s,x_1,u)-h^{ij}(s,x_2
            ,v)|+|\sigma(s,x_1,u)-\sigma(s,x_2,v)|\\
            &\leq C(|x_1-x_2|+|u-v|);\\
            &|b(s,x,u)|+|h^{ij}(s,x,u)|+|\sigma(s,x,u)|\leq C(1+|x|).
        \end{align*}
    \end{small}
\end{enumerate}
Then we turn to the assumptions on the coefficients of $G$-BSDE \eqref{state-2}.
\begin{enumerate}
    \item [{\bf (B1)}] $\Phi(x):\mathbb{R}^{n}\rightarrow\mathbb{R}$,
            $f(t,x,y,z,u),g^{ij}(t,x,y,z,u):[0,T]\times\mathbb{R}^{n}\times\mathbb{R}\times\mathbb{R}^{d}\times U\rightarrow\mathbb{R}$
            are deterministic functions, continuous in $t$ and uniformly
            continuous in $y,u$;
    \item [{\bf (B2)}] There is a constant $C>0$ such that for any 
            $t\in [0, T]$, $x_1,x_2\in\mathbb{R}^{n}$, $y\in\mathbb{R}$, $z_1,z_2\in\mathbb{R}^d,t\in[0,T]$
            $u\in U$
            \begin{small}
                \begin{align*}
                    &|f(t,x_1,y,z_1,u)-f(t,x_2,y,z_2,u)|
                    +|g^{ij}(t,x,y,z_1,u)-
                    g^{ij}(t,x,y,z_2,u)|\\
                    &+|\Phi(x_1)-\Phi(x_2)|\le C(|x_1-x_2|+|z_1-z_2|);
                \end{align*}
            \end{small}
    \item [{\bf (B3)}] There is a constant $\mu>0$
    such that for any $t\in [0, T]$, $x\in\mathbb{R}^{n}$,
    $y_1, y_2\in\mathbb{R}$, $z\in\mathbb{R}^d$, $u\in U$
    \begin{small}
        \begin{align*}
            (y_1-y_2)\Big[\Big(f(t,x,y_1,z,u)-f(t,x,y_2,z,u)\Big)+\Big(
            g^{ij}(t,x,y_1,z,u)-g^{ij}(t,x,y_2,z,u)\Big)\Big]
            \le \mu|y_1-y_2|^2;
        \end{align*}
    \end{small}
    \item [{\bf (B4)}] There is a constant $C>0$ such that for any 
    $t\in [0, T]$, $x\in\mathbb{R}^{n}$, $y\in\mathbb{R}$,
    $z\in\mathbb{R}^d$, $u\in U $
    \begin{small}
        \begin{align*}
            |f(t,x,y,z,u)-f(t,x,0,z,u)|+|g^{ij}(t,x,y,z,u)-g^{ij}(t,x,0,z,u)|\le C(1+|y|).
        \end{align*}
    \end{small}
\end{enumerate}
\begin{remark}
    Here, in addition to the monotonicity assumption, we weaken 
    the Lipschitz continuity of the aggregator 
    with respect to the control $u$ by the uniform continuity,
    which generalizes the results of Hu and Ji \cite{HJ17}. 
\end{remark}
We have the following theorems with the above assumptions.
\begin{theorem}
    (\cite{P10}) Let Assumptions {\bf (A1)} and {\bf (A2)} hold. 
    Then there exists a unique
    adapted solution $X$ for equation (\ref{state-1}).
\end{theorem}   
\begin{theorem}
    \label{estimate-state-1} (\cite{P10}) 
    Let 
    $\zeta,\zeta^{\prime}\in L_{G}^{p}(\Omega_{t};\mathbb{R}^{n})$ with $p\geq2$ and $u,v\in\mathcal{U}[t,T]$. 
    Then we
    have, for each $\delta\in[0,T-t]$,
    \begin{small}
       \[
            \widehat{\mathbb{E}}_{t}\Big[|X_{t+\delta}^{t,\zeta,u}-X_{t+\delta}^{t,\zeta^{\prime}
            ,v}|^{2}\Big]\leq C\Biggl(|\zeta-\zeta^{\prime}|^{2}+\widehat{\mathbb{E}}_{t}\biggl[\int
            _{t}^{t+\delta}|u_{s}-v_{s}|^{2}ds\biggr]\Biggr),
        \] 
    \end{small}
    \begin{small}
        \[
            \widehat{\mathbb{E}}_{t}\Big[|X_{t+\delta}^{t,\zeta,u}|^{p}\Big]\leq C\big(1+|\zeta|^{p}\big),
        \]
    \end{small}
    \begin{small}
        \[
            \widehat{\mathbb{E}}_{t}\Biggl[\sup_{s\in[ t,t+\delta]}|X_{s}^{t,\zeta,u}-\zeta
            |^{p}\Biggr]\leq C\big(1+|\zeta|^{p}\big)\delta^{p/2},
        \]
    \end{small}where the constant $C$ depends on $C$, $G$, $p$, $n$, $U$ and $T$.
\end{theorem}
For convenience, we suppose that $f$ and $g^{ij} $ are independent of $u$.
According to Lemma 1 in Lepeltier and Martin \cite{LSM97}, 
there exists a sequence of 
Lipschitz functions that nicely approximates 
$f$ and $g^{ij} $, respectively. 

For $\phi=f,g^{ij}$,
define 
\begin{small}
 \[
    \phi_n(t,x,y,z):=\inf_{q\in\mathbb{Q}}\{\phi(t,x,q,z)+n|y-q|\},
\]       
\end{small}then $\phi_n$ have the following properties.
\begin{lemma}\label{lemma4.1}
    Let $C$ be the constant in {\bf (B4)}. We have
    \begin{enumerate}
        
        \item [(i)] $\phi_n$ is well defined for $n\ge C$ and we have
        \begin{small}
        \[
            |\phi_n(t,x,y,z)|\le |\phi(t,0,0,0)|+C(1+|x|+|y|+|z|).
        \] 
        \end{small}   
        \item [(ii)] For any $x_1,x_2\in\mathbb{R}^n,y\in\mathbb{R},z_1,z_2\in\mathbb{R}^d,t\in[0,T],$ 
        \begin{small}
        \[
            |\phi_n(t,x_1,y,z_1)-\phi_n(t,x_2,y,z_2)|\le C(|x_1-x_2|+|z_1-z_2|).
        \]
        \end{small}
        \item [(iii)] For any $t\in [0, T],y_1, y_2\in\mathbb{R},z\in\mathbb{R}^d,$
        \begin{small}
        \[
            |\phi_n(t,x,y_1,z)-\phi_n(t,x,y_2,z)|\le n|y_1-y_2|.
        \]
        \end{small}
        \item [(iv)] For any $t\in [0, T],y\in\mathbb{R},z\in\mathbb{R}^d,$ the 
        sequence $(\phi_n(t,x,y,z))_n$ is increasing.
        \item [(v)] If $\phi$ is decreasing in $y$, then so is $\phi_n.$     
    \end{enumerate}
\end{lemma}

\begin{proof}
    The properties (i), (ii), (iii) and (iv) are either clear 
    or proved in \cite{LSM97}, (v) can be similarly proved in \cite{POD13}.
\end{proof}

Without loss of generality, 
we assume $g^{ij}\equiv 0$ for simplicity. Indeed, we have the following Lemma.

\begin{lemma}\label{lemma3.2}
    Let $\lambda>0,$ then $(Y_t, Z_t, K_t)$ solve the $G$-BSDE \eqref{state-2} if and only
    if $\small (e^{\lambda t}Y_t,e^{\lambda t}Z_t,\int_{0}^{t}e^{\lambda s}dK_s)$
    solve the $G$-BSDE with terminal condition $e^{\lambda T}\Phi(X_{T}^{t,\zeta})$ and 
    $f^\prime(t,x,y,z)=e^{\lambda t}f(t,x,e^{-\lambda t}y,e^{-\lambda t}z)-\lambda y$.
\end{lemma}

\begin{proof}
    It is easy to prove that the two solutions solve 
    the corresponding equations is a simple
    consequence of It\^{o}'s formula. Then by Lemma \ref{lemma2.6},  
    $\small \int_{0}^{t}e^{\lambda s}dK_s$ is a decreasing G-martingale. 
    And $\small \int_{0}^{T}e^{\lambda s}dK_s\in L_G^\alpha(\Omega_T)$ is obviously.
\end{proof}

If we choose $\lambda=\mu$ then $f^\prime$ satisfies the same 
assumptions as $f$, but with {\bf (B3)} replaced by
    \begin{enumerate}
        \item [{\bf ($\text{B3}^\prime$)}] $(y_1-y_2)[f(t,x,y_1,z)-f(t,x,y_2,z)]\le 0.$
    \end{enumerate}
Based on the above, we can assume without loss of generality that our generator
is decreasing in $y.$ Hence, we shall assume that $f$ satisfies {\bf (B1)}, {\bf (B2)}, {\bf ($\text{B3}^\prime$)}  and {\bf (B4)}.

We only need to slightly modify the Lemma 4.5 in \cite{Song2019}, then follow the simialr idea in \cite{Song2019} to obtain the following theorem. 
\begin{theorem}\label{our-well-posedness}
    Let assumptions {\bf (B1)}-{\bf (B4)} hold. Then there exists a unique adapted solution $(Y,Z,K)$ for equation (\ref{state-2}).
\end{theorem}
For a fixed $n,$ consider the following $G$-BSDE,
\begin{small}
    \begin{equation}
        Y_t^n=\Phi(X_{T})+\int_{t}^{T}f_n(s,X_s,Y_{s}^n,Z_{s}^n)ds
            -\int_{t}^{T}Z_{s}^ndB_{s}-(K_{T}^n-K_{t}^n).\label{state-n}
    \end{equation}
\end{small}Indeed, by Theorem \ref{estimate-state-1} and {\bf (B1)}-{\bf (B2)}, it is easy to obtain that $\Phi(X_{T})\in L_G^\beta(\Omega_T)$.
Moreover, according to Theorem 4.16 in \cite{HWZ16}, the fact  
$f(\cdot,X_{\cdot},y,z)\in M^\beta_G(0,T)$ is clear.
Therefore, by Lemma \ref{lemma4.1}
and Theorem \ref{thm2.7}, we konw that for all 
$n$ the above $G$-BSDE \eqref{state-n} has a
unique solution $(Y^n, Z^n, K^n) \in \mathfrak{S}_{G}^{\alpha}(0,T)$.
The idea to the proof of existence is to prove that 
the limit in a certain sense of the sequence
$(Y^n, Z^n, K^n)$ is a solution of the $G$-BSDE \eqref{state-2}.

Now, we prove the similar result as Lemma 4.5 in \cite{Song2019}.
\begin{proposition}
    For some $\beta>2$ and any $1<\alpha<\beta$, we have
    \begin{small}
        \begin{equation*}
            \sup_{0\leq t\leq T}\widehat{\mathbb{E}}\left[|Y_t^n-Y_t^m|^\alpha\right]\to 0, \text{ as } n,m\to\infty.
        \end{equation*}
    \end{small}
\end{proposition}
\begin{proof}
    Following the same derivation as Lemma 4.5 in \cite{Song2019} line by line, we 
    obtain that 
    \begin{align*}
        |Y_t^n-Y_t^m|^\alpha\leq C\int_{t}^{T}
                \widehat{\mathbb{E}}_t\Big[|\hat{f}_n|^{\alpha}+|\hat{f}_m|^{\alpha}\Big]ds,
    \end{align*}
    where 
    \begin{align*}
        \hat{f}_n=f(s,x,y,z)-f_n(s,x,y,z).
    \end{align*}
    Since $f$ and $f_n$ satisfy the linear growth condition, and by Theorem \ref{estimate-state-1} and 
    Proposition \ref{propA.1}-\ref{propA.2}, we have 
    \begin{small}
        \[
            \widehat{\mathbb{E}}\Big[|\hat{f}_n|^{2+\epsilon}\Big]\leq C, 
        \]
    \end{small}for some $C$ independent of $n$ and any $\epsilon>0$.
    Moreover, we can check that $\hat{f}_n$ converges to $0$.

    Then, it is easy to check that $X\in L_G^1(\Omega_T)$ ($X$ is the solution of $G$-SDE \eqref{state-1}),
    then by Theorem \ref{thm-LGp}, $X$ has a quasi-continuous version. Recalling 
    $f(\cdot,X_{\cdot},y,z)\in M^\beta_G(0,T)$, by 
    the characterization of $M^\beta_G(0,T)$, i.e. Theorem \ref{thm-MGp}, we can obtain that 
    $f(\cdot,X_{\cdot},y,z)$ has a quasi-continuous version. Thus, $f_n(\cdot,X_{\cdot},y,z)$ also has a quasi-continuous version.
    These imply that $\hat{f}_n$ has a quasi-continuous version.

    Finally, applying Markov inequality and H\"{o}lder inequality under $G$-framework,
    we get 
    \begin{small}
        \begin{align*}
            \lim_{N\to\infty}\widehat{\mathbb{E}}\Big[|\hat{f}_n|^\alpha\mathbf{1}_{|\hat{f}_n|^\alpha>N}\Big]=0.
        \end{align*}
    \end{small}Based the above results and 
    by Theorem \ref{thm-LGp},
    we have $|\hat{f}_n|^\beta\in L_G^1(\Omega_T)$.
    Using the monotone convergence Theorem \ref{monotone-con-thm}, we have 
    \begin{small}
        \begin{equation*}
            \widehat{\mathbb{E}}\Big[|\hat{f}_n|^\alpha\Big]\to 0, \text{ as } n\to +\infty, \text{ for all } 0\leq t\leq T.
        \end{equation*} 
    \end{small}
    Similar analysis holds for $\small \widehat{\mathbb{E}}\Big[|\hat{f}_m|^\alpha\Big]$.    
\end{proof}

Next, we can borrow the same idea as Theorem 4.1-4.2 in \cite{Song2019} to obtain the Theorem \ref{our-well-posedness}.

By Proposition \ref{propA.3} and Remark \ref{remarkA.4}, we have the following
estimate, which differs from the Proposition 5.1 in \cite{HJPS1}
only in the estimated coefficients. 

\begin{theorem}\label{estimate-state-2}
    Let 
    $\zeta,\zeta^{\prime}\in\cup_{\varepsilon>0}L_{G}^{2+\varepsilon}(\Omega_{t};\mathbb{R}^{n})$ 
    and $u,v\in\mathcal{U}[t,T]$. 
    Then there exist a positive constant $C$ depending 
    on $C$, $G$, $\mu$ and $T$ such that
    \begin{small}
        \begin{align*}
            |Y_{t}^{t,\zeta,u}-Y_{t}^{t,\zeta^{\prime},v}|^{2}  
            \leq C\widehat{\mathbb{E}}_{t}\Biggl[|\Phi(X_{T}^{t,\zeta,u})-\Phi
            (X_{T}^{t,\zeta^{\prime},v})|^{2}+\int_{t}^{T}|\hat{F}_{s}|^{2}ds\Biggr],
        \end{align*}
    \end{small}where
    \begin{small}
        \begin{align*}
            \hat{F}_{s}  &  =|f(s,X_{s}^{t,\zeta,u},Y_{s}^{t,\zeta,u},Z_{s}^{t,\zeta,u}
            ,u_{s})-f(s,X_{s}^{t,\zeta^{\prime},v},Y_{s}^{t,\zeta,u},Z_{s}^{t,\zeta,u},v_{s})|\\
            &  +\sum_{i,j=1}^{d}|g^{ij}(s,X_{s}^{t,\zeta,u},Y_{s}^{t,\zeta,u},Z_{s}^{t,\zeta
            ,u},u_{s})-g^{ij}(s,X_{s}^{t,\zeta^{\prime},v},Y_{s}^{t,\zeta,u},Z_{s}^{t,\zeta
            ,u},v_{s})|.
        \end{align*}
    \end{small} 
\end{theorem}
In the end of this subsection, we provide the our comparison
theorem which is crucial to obtain the dynamic programming principle.
Due to the non-Lipschitz assumptions, it is different from 
the comparison Theorem 3.6 in \cite{HJb}, we can not adopt the linearization 
method. We refer to the idea of Theorem 3.7 in \cite{He22} to prove 
our comparison theorem. For convenience, we assume that the generators are independent of $x$ and $u$.
\begin{theorem}[Comparison Theorem]\label{thm3.9}
    Suppose $\xi^l\in L_{G}^{\beta}(\Omega_T),l=1,2,$ for some $\beta>2$, and 
    $f^l,g^{l,ij}$ satisfy assumptions {\bf (B1)}, {\bf (B2)}, {\bf ($\text{B3}^\prime$)} and {\bf (B4)}.
    Let $(Y^l, Z^l, K^l)\in \mathfrak{S}_{G}^{2 }(0,T)$ 
    be a solution of the G-BSDE \eqref{state-2} corresponding to the data
    $(\xi^l,f^l,g^{l,ij})$. If $\xi^1\le\xi^2$, $f^1(t,y,z)\le f^2(t,y,z)$, 
    $g^{1,ij}(t,y,z)\le g^{2,ij}(t,y,z)$ for any $(t,\omega,y,z)$, then 
    $Y_t^1\le Y_t^2$ for each $t$.
\end{theorem}
\begin{proof}
    Without loss of generality, we assume $g^{ij} \equiv 0$.
    Set $\hat{Y}_t=Y^1_t-Y^2_t$, 
    $\hat{Z}_t=Z^1_t-Z^2_t$, $\hat{K}_t=K_t^1-K_t^2$,
    $\hat{\xi}=\xi^1-\xi^2\le 0$,
    $\hat{f}_t=f^1(s,Y_t^1,Z^1_t)-f^2(s,Y_t^2,Z_t^2)$, we have 
    \begin{small}
        \begin{equation*}
            \hat{Y}_t=\hat{\xi}+\int_{t}^{T}\hat{f}_sds-\int_{t}^{T}
            \hat{Z}_sdB_s-\int_{t}^{T}d\hat{K}_s.
        \end{equation*}
    \end{small}By the representation Theorem \ref{thm2.1} of $G$-expectation,
    there exists a weakly compact family of probability measures
    $\mathcal{P}$. For each $P\in\mathcal{P}$, $(\hat{Y}_t,\hat{Z}_t)$
    solves the following classical BSDE:
    \begin{small}
        \begin{equation*}
            \hat{Y}_t=\hat{\xi}-(\hat{K}_T-\hat{K}_t)+\int_{t}^{T}\hat{f}_sds-\int_{t}^{T}
            \hat{Z}_sdB_s,\ P-a.s..
        \end{equation*}
    \end{small}Note that $\hat{Y}_t$ is a continuous semimartingale.
    According to the Tanaka-Meyer formula, we get
    \begin{small}
        \begin{equation*}
            \hat{Y}_t^+=\hat{\xi}^+-\int_{t}^{T}\mathbf{1}_{\{\hat{Y}_s>0\}}d\hat{Y}_s-\frac{1}{2}\int_{t}^{T}dL_s^P(0),\ P-a.s.,
        \end{equation*}
    \end{small}where $L^P_t (0)$ is the local time of $\hat{Y}_t$ at $0$ under measure $P$.
    Applying It\^{o}'s formula to $|\hat{Y}_t^+|^2$, we get
    \begin{small}
        \begin{align*}
            &|\hat{Y}_t^+|^2+\int_{t}^{T}\mathbf{1}_{\{\hat{Y}_s>0\}}|\hat{Z}_s|^2
            d\langle B\rangle_s+\int_{t}^{T}\hat{Y}_s^+dL_s^P(0)\\
            &=|\hat{\xi}^+|^2+2\int_{t}^{T}\hat{Y}_s^+\hat{f}_sds-2\int_{t}^{T}
            \Big(\hat{Y}_s^+\hat{Z}_sdB_s+\hat{Y}_s^+d(K_s^1-K_s^2)\Big), \ P-a.s..
        \end{align*}
    \end{small}Moreover, $\small \int_{0}^{t}\mathbf{1}_{\{\hat{Y}_s>0\}}dL_s^P(0)=0$ (see \cite{KS14}),
    thus $\small \int_{t}^{T}\hat{Y}_s^+dL_s^P(0)=0$.
    Note that $\hat{\xi}^+=0$, then we obtain
    \begin{small}
        \begin{align*}
            |\hat{Y}_t^+|^2+\int_{t}^{T}\mathbf{1}_{\{\hat{Y}_s>0\}}|\hat{Z}_s|^2
            d\langle B\rangle_s+2(J_T-J_t)\leq 
            2\int_{t}^{T}\hat{Y}_s^+\mathbf{1}_{\{\hat{Y}_s>0\}}\hat{f}_sds,\ P-a.s.,
        \end{align*}
    \end{small}where $\small J_t=\int_{0}^{t}\hat{Y}_s^+\hat{Z}_sdB_s+\int_{0}^{t}\hat{Y}_s^+dK_s^1$
    is a $G$-martingale.

    Using the monotonic decreasing 
    property of $f$ in $y$
    and $f^1\leq f^2$, we have
    \begin{small}
        \begin{align*}
            2\int_{t}^{T}\hat{Y}_s^+\mathbf{1}_{\{\hat{Y}_s>0\}}\hat{f}_sds
            \leq \int_{t}^{T}2C|\hat{Y}_s^+|\mathbf{1}_{\{\hat{Y}_s>0\}}|\hat{Z}_s|ds
            \leq \frac{C^2}{\underline{\sigma}^2}\int_{t}^{T}|\hat{Y}_s^+|^2ds+
            \int_{t}^{T}\mathbf{1}_{\{\hat{Y}_s>0\}}|\hat{Z}_s|^2d\langle B\rangle_s.
        \end{align*}
    \end{small}Then it follows that 
    \begin{small}
        \begin{align*}
            |\hat{Y}_t^+|^2+2(J_T-J_t)\leq 
            \frac{C^2}{\underline{\sigma}^2}\int_{t}^{T}|\hat{Y}_s^+|^2ds,\ P-a.s..
        \end{align*}
    \end{small}Taking $G$-conditional expectation first, then taking $G$-expectation, we obtain
    \begin{small}
        \begin{equation*}
            \widehat{\mathbb{E}}\Big[|\hat{Y}_t^+|^2\Big]\leq 
            \frac{C^2}{\underline{\sigma}^2}\int_{t}^{T}\widehat{\mathbb{E}}\Big[|\hat{Y}_s^+|^2\Big]ds.
        \end{equation*}
    \end{small}By Gronwall inequality, we get 
    \begin{small}
        \[
            \widehat{\mathbb{E}}\Big[|\hat{Y}_t^+|^2\Big]\le 0.
        \]
    \end{small}It follows that 
    \begin{small}
        \[
            \widehat{\mathbb{E}}\Big[|\hat{Y}_t^+|^2\Big]= 0.
        \] 
    \end{small}Hence, by the continuity of $Y$, we deduce that $Y_t^1\le Y_t^2$ for each $t$.
\end{proof}

\subsection{Problem formulation}
In this subsection, we mainly give the specific formulation
of the stochastic recursive optimal control problem driven by $G$-Brownian motion.
The state equation is governed by $G$-SDE \eqref{state-1}, and the 
objective functional is described by the solution of the $G$-BSDE \eqref{state-2}
at time $t$. Our stochastic recursive optimal control problem is 
to find $u\in \mathcal{U}[t,T]$ which minimizes the 
objective function $Y_t^{t,x,u}$ for each given 
$x\in\mathbb{R}^n$. To this end, we need to define the
essential infimum of $\big\{Y_{t}^{t,x,u}\mid u\in\mathcal{U}[t,T]\big\}.$
\begin{definition}\label{essinf}
    The essential infimum of 
    $\{Y_{t}^{t,x,u}\mid u\in\mathcal{U}[t,T]\}$, denoted by 
    $\small {\text{ess}\inf}_{u(\cdot)\in\mathcal{U}[t,T]}\ Y_{t}^{t,x,u}$, 
    is a random variable $\eta\in L_{G}^{2}(\Omega_{t})$ satisfying:
    \begin{enumerate}
        \item [(i)]$\forall u\in\mathcal{U}[t,T],$ $\eta\leq Y_{t}^{t,x,u}$ q.s.;
        \item [(ii)]if $\eta^\prime$ is a random variable satisfying 
        $\eta^\prime\leq Y_{t}^{t,x,u}$
        q.s. for any $u\in\mathcal{U}[t,T]$, then $\eta\geq \eta^\prime$ q.s..
    \end{enumerate}
\end{definition}
For $\zeta\in\cup_{\varepsilon>0}L_{G}^{2+\varepsilon}(\Omega_{t};\mathbb{R}^{n})$,
we can define the essential infimum of 
$\big\{Y_{t}^{t,\zeta,u}\mid u\in\mathcal{U}[t,T]\big\}$ 
by similar way.
\begin{remark}\label{remark-essinf-not-exist}
    The essential infimum may be not exist, see Example 11 in \cite{HJ17}.
\end{remark}
Finally, for each $(t,x)\in[0,T]\times\mathbb{R}^{n}$, 
the value function our problem is defined by
\begin{small}
    \begin{equation}\label{valuefunction0}
        V(t,x):=\underset{u\in\mathcal{U}[t,T]}{ess\inf}\ Y_{t}^{t,x,u}. 
    \end{equation}
\end{small}
\subsection{Dynamic programming principle}
As mentioned in Remark \ref{remark-essinf-not-exist}, 
to obtain our dynamic programming principle, we first prove 
some properties of the value function, i.e.,
the existence and determinacy of $V(\cdot,\cdot)$, and for each 
$\zeta\in\cup_{\varepsilon>0}L_{G}^{2+\varepsilon}(\Omega_{t};\mathbb{R}^{n})$, 
$\small V(t,\zeta)={ess\inf}_{u\in\mathcal{U}[t,T]}\ Y_{t}^{t,\zeta,u}$.

We use the following notations:
\begin{small}
    \begin{align*}
        L_{ip}(\Omega_{s}^{t})  &  :=\big\{ \varphi(B_{t_{1}}-B_{t},...,B_{t_{n}}
        -B_{t}):n\geq1,t_{1},...,t_{n}\in[t,s],\varphi\in C_{b.Lip}
        (\mathbb{R}^{d\times n})\big\};\\
        L_{G}^{2}(\Omega_{s}^{t})  &  :=\Big\{ \text{the completion of }L_{ip}(\Omega
        _{s}^{t})\text{ under the norm }\Vert\cdot\Vert_{2,G}\Big\};\\
        M_{G}^{0,t}(t,T)  &  :=\Biggl\{ \eta_{s}=\sum_{i=0}^{N-1}\xi_{i}I_{[t_{i},t_{i+1}
        )}(s):t=t_{0}<\cdots<t_{N}=T,\xi_{i}\in L_{ip}(\Omega_{t_{i}}^{t})\Biggr\};\\
        M_{G}^{2,t}(t,T)  &  :=\Big\{ \text{the completion of }M_{G}^{0,t}(t,T)\text{
        under the norm }\Vert\cdot\Vert_{M_{G}^{2}}\Big\};\\
        \mathcal{U}[t,T]  &  :=\big\{u:u\in M_{G}^{2}(t,T;\mathbb{R}^{m})\text{ with
        values in }U\big\};\\
        \mathcal{U}^{t}[t,T]  &  :=\big\{u:u\in M_{G}^{2,t}(t,T;\mathbb{R}^{m})\text{ with
        values in }U\big\};\\
        \mathbb{U}[t,T]  &  :=\Biggl\{u=\sum\limits_{i=1}^{n}I_{A_{i}}u^{i}:n\in
        \mathbb{N},u^{i}\in\mathcal{U}^{t}[t,T],I_{A_{i}}\in L_{G}^{2}(\Omega
        _{t}),\Omega=
        {\displaystyle\bigcup\limits_{i=1}^{n}}A_{i}\Biggr\};\\
        \mathbb{U}^{t}[t,T]  &  :=\Biggl\{u=\sum_{i=0}^{N-1}(\sum_{j=1}^{l_{i}}a_{j}
        ^{i}I_{A_{j}^{i}})I_{[t_{i},t_{i+1})}(s):l_{i}\in\mathbb{N},a_{j}^{i}\in
        U,I_{A_{j}^{i}}\in L_{G}^{2}(\Omega_{t_{i}}^{t}),
        \Omega={\displaystyle\bigcup\limits_{j=1}^{l_{i}}}
        A_{j}^{i}\Biggr\}.
    \end{align*}
\end{small}
\begin{remark}
    Note that 
    $\mathbb{U}^{t}[t,T]\subset\mathcal{U}^{t}[t,T]\subset\mathbb{U}[t,T]\subset\mathcal{U}[t,T]$.
\end{remark}
In order to show $V(\cdot,\cdot)$ is deterministic, we recall the following lemmas.
\begin{lemma}[\cite{HJ17}]
    \label{u-uk-1} Let $u\in\mathcal{U}[t,T]$ be given. Then there exists a
    sequence $(u^{k})_{k\geq1}$ in $\mathbb{U}[t,T]$ such that
    \begin{small}
    \[
        \lim_{k\rightarrow\infty}\widehat{\mathbb{E}}\Biggl[\int_{t}^{T}|u_{s}-u_{s}^{k}%
        |^{2}ds\Biggr]=0.
    \]
    \end{small}
\end{lemma}

\begin{lemma}[\cite{HJ17}]
    \label{u-uk-2}Let $u\in\mathcal{U}^{t}[t,T]$ be given. 
    Then there exists a
    sequence $(u^{k})_{k\geq1}$ in $\mathbb{U}^{t}[t,T]$ such that
    \begin{small}
    \[
    \lim_{k\rightarrow\infty}\widehat{\mathbb{E}}\Biggl[\int_{t}^{T}|u_{s}-u_{s}^{k}
    |^{2}ds\Biggr]=0.
    \]
    \end{small}
\end{lemma}

\begin{theorem}\label{deterministic-V} 
    The value function $V(t,x)$ exists and
    \begin{small}
    \[
        V(t,x)=\inf_{u\in\mathcal{U}^{t}[t,T]}Y_{t}^{t,x,u}=\inf_{u\in\mathbb{U}
        ^{t}[t,T]}Y_{t}^{t,x,u}.
    \]
    \end{small}
\end{theorem}

\begin{proof}
    For each $v\in\mathcal{U}^{t}[t,T]$, 
    we can easily check that $Y_{t}^{t,x,v}$
    is a constant. In the following, we prove that 
    $\small {ess\inf}_{u\in\mathcal{U}[t,T]}\ Y_{t}^{t,x,u}=\inf_{v\in\mathcal{U}^{t}[t,T]}Y_{t}^{t,x,v},$
    q.s.. 
    Note that $\mathcal{U}^{t}[t,T]\subset\mathcal{U}[t,T]$, 
    we only need to
    verify that $\small Y_{t}^{t,x,u}\geq\inf_{v\in\mathcal{U}^{t}[t,T]}Y_{t}^{t,x,v}$, q.s.
    for each $u\in\mathcal{U}[t,T]$.

    Then by Lemma \ref{u-uk-1}, for each $u\in\mathcal{U}[t,T]$, 
    there exists a sequence 
    $\small u^k=\sum_{i=1}^{N_k}I_{A_{i}^k}v^{i,k}\in\mathbb{U}[t,T]$, 
    such that $\tiny \widehat{\mathbb{E}}\Biggl[\int_{t}^{T}\biggl|u_{s}-u_s^k\biggr|^{2}ds\Biggr]\rightarrow0.$
    We consider the following equations:
    \begin{small}
        \begin{align*}
            \left\{
            \begin{array}{l}
                dX_{s}^{t,\zeta,v^{i,k}}=b(s,X_{s}^{t,\zeta,v^{i,k}},v_{s}^{i,k})ds+h^{ij}
                (s,X_{s}^{t,\zeta,v^{i,k}},v_{s}^{i,k})d\langle B^{i},B^{j}\rangle_{s}
                +\sigma(s,X_{s}^{t,\zeta,v^{i,k}},v_{s}^{i,k})dB_{s},\\
                X_{t}^{t,\zeta,v^{i,k}}=\zeta,\ s\in [t,T]\text{, }i=1,...,N.           
            \end{array}\right.
        \end{align*}
    \end{small}
    \begin{small}
        \begin{align*}
            Y_{s}^{t,\zeta,v^{i,k}}  &  =\Phi(X_{T}^{t,\zeta,v^{i,k}})+\int_{s}^{T}f(r,X_{r}^{t,\zeta,v^{i,k}},Y_{r}^{t,\zeta,v^{i,k}},Z_{r}^{t,\zeta
            ,v^{i,k}},v^{i,k}_{r})dr\notag\\
            &\quad+\int_{s}^{T}g^{ij}(r,X_{r}^{t,\zeta,v^{i,k}},Y_{r}^{t,\zeta,v^{i,k}},Z_{r}^{t,\zeta
            ,v^{i,k}},v^{i,k}_{r})d\langle B^{i},B^{j}\rangle_{r}\notag\\
            &\quad  -\int_{s}^{T}Z_{r}^{t,\zeta,v^{i,k}}dB_{r}-(K_{T}^{t,\zeta,v^{i,k}}-K_{s}^{t,\zeta,v^{i,k}}),\ \ s\in[
                t,T]\text{, }i=1,...,N.
        \end{align*}
    \end{small}For $s\in [t,T]$, we set
    \begin{small}
    \[
        \bar{X}_{s}^{t,\xi,u^k}=\sum_{i=1}^{N_k}I_{A_{i}^k}X_{s}^{t,\xi,v^{i,k}},\ \bar{Y}
        _{s}^{t,\xi,u^k}=\sum_{i=1}^{N_k}I_{A_{i}^k}Y_{s}^{t,\xi,v^{i,k}},\ \bar{Z}_{s}
        ^{t,\xi,u^k}=\sum_{i=1}^{N_k}I_{A_{i}^k}Z_{s}^{t,\xi,v^{i,k}},\ \bar{K}_{s}^{t,\xi
        ,u^k}=\sum_{i=1}^{N_k}I_{A_{i}^k}K_{s}^{t,\xi,v^{i,k}}.
    \]
    \end{small}Multiplying $I_{A_{i}}$ on both sides of the above two equations and summing up,
    we have
    \begin{small}
        \begin{equation*}
            \left\{
            \begin{array}{l}
                d\bar{X}_{s}^{t,\zeta,u^k}=b(s,\bar{X}_{s}^{t,\zeta,v},u^k_{s})ds+h^{ij}
                (s,\bar{X}_{s}^{t,\zeta,u^k},u^k_{s})d\langle B^{i},B^{j}\rangle_{s}
                +\sigma(s,\bar{X}_{s}^{t,\zeta,u^k},u^k_{s})dB_{s},\\
                X_{t}^{t,\zeta,u^k}=\zeta,\ s\in [t,T],          
            \end{array}\right.
        \end{equation*}
    \end{small}
    \begin{small}
        \begin{align*}
            \bar{Y}_{s}^{t,\zeta,u^k}  &  =\Phi(\bar{X}_{T}^{t,\zeta,u^k})+\int_{s}^{T}f(r,\bar{X}_{r}^{t,\zeta,u^k},\bar{Y}_{r}^{t,\zeta,u^k},\bar{Z}_{r}^{t,\zeta
            ,u^k},u^k_{r})dr\notag\\
            &\quad+\int_{s}^{T}g^{ij}(r,\bar{X}_{r}^{t,\zeta,u^k},\bar{Y}_{r}^{t,\zeta,u^k},\bar{Z}_{r}^{t,\zeta
            ,u^k},u^k_{r})d\langle B^{i},B^{j}\rangle_{r}\notag\\
            &\quad  -\int_{s}^{T}\bar{Z}_{r}^{t,\zeta,u^k}dB_{r}-(\bar{K}_{T}^{t,\zeta,u^k}-\bar{K}_{s}^{t,\zeta,u^k}),\ \ s\in[
                t,T].
        \end{align*}
    \end{small}By Theorem \ref{estimate-state-2} and Theorem \ref{estimate-state-1}, we get 
\begin{small}
    \begin{align*}
        &\widehat{\mathbb{E}}\Biggl[\biggl|Y_{t}^{t,\zeta,u}-\sum_{i=1}^{N_k}I_{A_{i}^k}Y_{t}^{t,\zeta,v^{i,k}
            }\biggr|^{2}\Biggr]\notag\\
        &\leq C \Biggl\{\widehat{\mathbb{E}}\Biggl[\int_{t}^{T}|u_s-u^k_s|^2ds\Biggr]
        +\widehat{\mathbb{E}}\Biggl[\int_{t}^{T}\Big|f(s,X_{s}^{t,\zeta,u},Y_{s}^{t,\zeta,u},Z_{s}^{t,\zeta,u}
        ,u_{s})-f(s,X_{s}^{t,\zeta,u},Y_{s}^{t,\zeta,u},Z_{s}^{t,\zeta,u},u^k_{s})\Big|^2ds\Biggr]\notag\\
        &\qquad +\sum_{i,j=1}^{d}\widehat{\mathbb{E}}\Biggl[\int_{t}^{T}\Big|g^{ij}(s,X_{s}^{t,\zeta,u},Y_{s}^{t,\zeta,u},Z_{s}^{t,\zeta
        ,u},u_{s})-g^{ij}(s,X_{s}^{t,\zeta,u},Y_{s}^{t,\zeta,u},Z_{s}^{t,\zeta
        ,u},u^k_{s})\Big|^2ds\Biggr]\Biggr\}.\label{Yu-barY}
    \end{align*}
\end{small}We only give the analysis for the second term of the RHS of the above inequality,
the third term can be treated similarly.
By the uniform continuity of $f$ in $u$, for any fixed $\epsilon> 0$, 
there exists a constant $\delta > 0$, such that 
$|f(\cdot,\cdot,\cdot,\cdot,u_1)-f(\cdot,\cdot,\cdot,\cdot,u_2)|<\epsilon$
whenever $|u_1-u_2| \le \delta$. Then we obtain that
\begin{small}
    \begin{align*}
        &\widehat{\mathbb{E}}\Biggl[\int_{t}^{T}\Big|f(s,X_{s}^{t,\zeta,u},Y_{s}^{t,\zeta,u},Z_{s}^{t,\zeta,u}
        ,u_{s})-f(s,X_{s}^{t,\zeta,u},Y_{s}^{t,\zeta,u},Z_{s}^{t,\zeta,u},u^k_{s})\Big|^2ds\Biggr]\\
        &\leq 2 \widehat{\mathbb{E}}\Biggl[\int_{t}^{T}\Big|f(s,X_{s}^{t,\zeta,u},Y_{s}^{t,\zeta,u},Z_{s}^{t,\zeta,u}
        ,u_{s})-f(s,X_{s}^{t,\zeta,u},Y_{s}^{t,\zeta,u},Z_{s}^{t,\zeta,u},u^k_{s})\Big|^2\mathbf{1}_{\{|u_s-u^k_s|\le\delta\}}ds\Biggr]\\
        &\quad+2 \widehat{\mathbb{E}}\Biggl[\int_{t}^{T}\Big|f(s,X_{s}^{t,\zeta,u},Y_{s}^{t,\zeta,u},Z_{s}^{t,\zeta,u}
        ,u_{s})-f(s,X_{s}^{t,\zeta,u},Y_{s}^{t,\zeta,u},Z_{s}^{t,\zeta,u},u^k_{s})\Big|^2\mathbf{1}_{\{|u_s-u^k_s|>\delta\}}ds\Biggr]\\
        &\le 2(T-t)\epsilon^2+2\widehat{\mathbb{E}}\Biggl[\int_{t}^{T}\biggl(2C(1+|Y_{s}^{t,\zeta,u}|)+2(|X_{s}^{t,\zeta,u}|+|Z_{s}^{t,\zeta,u}|)\\
        &\quad+|f(s,0,0,0,u_s)|+|f(s,0,0,0,u^k_s)|\biggr)^2\mathbf{1}_{\{|u_s-u^k_s|>\delta\}}ds\Biggr],
    \end{align*}
\end{small}Note that $u^k$ also converges to $u$ in $M_{G}^{2}(0,T)$, 
    we have $\small \lim_{k\to\infty} \tiny \widehat{\mathbb{E}}\biggl[\int_{0}^{T}\mathbf{1}_{\{|u_s-u^k_s|> \delta\}}ds\biggr]=0,$
    and it is easy to verify that 
    $2C(1+|Y_{s}^{t,\zeta,u}|)+2(|X_{s}^{t,\zeta,u}|+|Z_{s}^{t,\zeta,u}|)+|f(s,0,0,0,u_s)|+|f(s,0,0,0,u^k_s)|\in M_{G}^{2}(0,T)$. 
    Then similar as the {\bf Step 3} in Theorem 3.7 of \cite{WZ21},
\begin{small}
    \begin{align*}
        &\lim_{k\to\infty}2\widehat{\mathbb{E}}\Biggl[\int_{t}^{T}\biggl(2C(1+|Y_{s}^{t,\zeta,u}|)+2(|X_{s}^{t,\zeta,u}|+|Z_{s}^{t,\zeta,u}|)\\
        &\qquad\qquad\quad+|f(s,0,0,0,u_s)|+|f(s,0,0,0,u^k_s)|\biggr)^2\mathbf{1}_{\{|u_s-u^k_s|>\delta\}}ds\Biggr]=0.
    \end{align*}
\end{small} 

Hence, we have
\begin{small}
    \[
    \lim_{k\to\infty}\widehat{\mathbb{E}}\Biggl[\biggl|Y_{t}^{t,\zeta,u}-\sum_{i=1}^{N_k}I_{A_{i}^k}
    Y_{t}^{t,x,v^{i,k}}\biggr|^{2}\Biggr]=0,
    \]
\end{small}which implies that $\small \sum_{i=1}^{N_k}I_{A_{i}^k}Y_{t}^{t,x,v^{i,k}}$ converges
to $Y_{t}^{t,x,u}$ in $\mathbb{L}_{G}^{2}$.

Then there exists a subsequence,
we still denote it by 
$\small \big\{ \sum_{i=1}^{N_{k}}I_{A_{i}^{k}}Y_{t}^{t,x,v^{i,k}}\big\}$,
which converges to $Y_{t}^{t,x,u}$ q.s.. 
We observed that
\begin{small}
    \[
    \sum_{i=1}^{N_{k}}I_{A_{i}^{k}}Y_{t}^{t,x,v^{i,k}}\geq\inf_{v\in
    \mathcal{U}^{t}[t,T]}Y_{t}^{t,x,v}\text{ q.s.},
    \]
\end{small}then $\small Y_{t}^{t,x,u}\geq\inf_{v\in\mathcal{U}^{t}[t,T]}Y_{t}^{t,x,v}$
q.s.. Thus
\begin{small}
    \[
    V(t,x)=\underset{u\in\mathcal{U}[t,T]}{ess\inf}\ Y_{t}^{t,x,u}=\inf
    _{v\in\mathcal{U}^{t}[t,T]}Y_{t}^{t,x,v}.
    \]
\end{small}By Lemma \ref{u-uk-2} and similar analysis mentioned above, 
we can conclude that
\begin{small}
    \[
    \inf_{u\in\mathcal{U}^{t}[t,T]}Y_{t}^{t,x,u}=\inf_{u\in\mathbb{U}^{t}[t,T]}Y_{t}^{t,x,u}.
\]
\end{small}Consequently, we finish the proof.
\end{proof}

The following lemma indicates that the value function $V(t, x)$ 
is Lipschitz continuous and linear growth with respect to $x$.
\begin{lemma}\label{Vx-Vy}
    There exists a constant $C>0$ depending on $T$, $G$ and
    $C$ such that for any $x,y\in\mathbb{R}^{n}$
    \begin{enumerate}
        \item [(i)] $\mid V(t,x)-V(t,y)\mid\leq C\mid x-y\mid$;
        \item [(ii)] $|V(t,x)|\le C(1+|x|)$.
    \end{enumerate}
\end{lemma}
\begin{proof}
    By Theorem \ref{estimate-state-2}, 
    for for any $x$, $y\in\mathbb{R}^{n}$ and $u\in\mathcal{U}^{t}[t,T]$,
    we have
    \begin{small}
        \begin{align*} 
            \Big| Y_{t}^{t,x,u}-Y_{t}^{t,y,u}\Big|^{2}
            \leq  C\widehat{\mathbb{E}}_t\Biggl[\Big|X_{T}^{t,x,u}-X_{T}^{t,y,u}
            \Big|^{2}+\int_{t}^{T}\Big| X_{s}^{t,x,u}-X_{s}^{t,y,u}\Big|^{2}ds\Biggr],   
    \end{align*} 
    \end{small}then by Theorem \ref{estimate-state-1}, we get 
    \begin{small}
    \[
        \widehat{\mathbb{E}}_t\biggl[\Big| X_{s}^{t,x,u}-X_{s}^{t,y,u}\Big|^{2}\biggr]
        \leq C|x-y|^2,\ s\in[t,T].
    \]
    \end{small}Furthermore, we have 
    \begin{small}
    \[
        \Biggl|\inf_{v\in\mathcal{U}^{t}[t,T]}Y_{t}^{t,x,v}
        -\inf_{v\in\mathcal{U}^{t}[t,T]}Y_{t}^{t,y,v}\Biggr|\leq\sup_{v\in\mathcal{U}
        ^{t}[t,T]}\Big|Y_{t}^{t,x,v}-Y_{t}^{t,y,v}\Big|.
    \]
    \end{small}Thus, by the above three inequalities, we prove the result (i).
    Similarly, we can obtain (ii).
\end{proof}

The following theorem will be a key role in the 
derivation of dynamic programming principle.
\begin{theorem}\label{V-zeta}
    For any 
    $\zeta\in\cup_{\varepsilon>0}L_{G}^{2+\varepsilon}(\Omega_{t};\mathbb{R}^{n})$, we have
    \begin{small}
    \[
        V(t,\zeta)=\underset{u\in\mathcal{U}[t,T]}{ess\inf}\ Y_{t}^{t,\zeta,u}.
    \]
    \end{small}
\end{theorem}
\begin{proof}
    The proof is similar to the one in \cite{HJ17}.
    For convenience, we provide the sketch of the proof.

    We first show that $\forall u\in\mathcal{U}[t,T]$, 
    $V(t,\zeta)\leq Y_{t}^{t,\zeta,u}\;$ q.s..

    For a fixed 
    $\zeta\in\cup_{\varepsilon>0}L_{G}^{2+\varepsilon}(\Omega_{t};\mathbb{R}^{n})$, 
    there exists a sequence $\small \zeta^{k}=\sum_{i=1}^{N_{k}}x_{i}^{k}I_{A_{i}^{k}}$, $k=1,2,...,$
    where $x_{i}^{k}\in\mathbb{R}^{n}$ and $\{A_{i}^{k}\}_{i=1}^{N_{k}}$ 
    is a $\mathcal{B}(\Omega_{t})$-partition of $\Omega$, such that
    \begin{small}
        \begin{equation}\label{zeta-zetak}
            \lim_{k\rightarrow\infty}\widehat{\mathbb{E}}\Big[|\zeta-\zeta^{k}|^{2}\Big]=0.   
        \end{equation}
    \end{small}By Lemma \ref{Vx-Vy}, we get
    \begin{small}
        \begin{equation}\label{Vzeta-Vzetak}
            \lim_{k\rightarrow\infty}\widehat{\mathbb{E}}\Big[\big|V(t,\zeta)-V(t,\zeta^{k})\big|^{2}\Big]=0.
        \end{equation}
    \end{small}
    Applying similar method as in the proof of Theorem \ref{deterministic-V}, 
    then by Theorem \ref{estimate-state-2}, we have 
    \begin{small}
        \begin{align}\label{Y-zeta-Y-zetak}
            & \Biggl|Y_{t}^{t,\zeta,u}-\sum_{i=1}^{N_{k}}I_{A_{i}^{k}}Y_{t}^{t,x_{i}^{k},u}\Biggr|^{2}\notag\\
            \leq & C\widehat{\mathbb{E}}_t\Biggl[\Biggl|\Phi(X_{T}^{t,\zeta,u})-\Phi\biggl(\sum
            _{i=1}^{N_{k}}I_{A_{i}^{k}}X_{T}^{t,x_{i}^{k},u}\biggr)\Biggr|^{2}+\int_{t}
            ^{T}\biggl| X_{s}^{t,\zeta,u}-\sum_{i=1}^{N_{k}}I_{A_{i}^{k}}X_{s}^{t,x_{i}^{k}
            ,u}\biggr|^{2}ds\Biggr]\notag\\
            \leq & C\widehat{\mathbb{E}}_t\Big[\big|\zeta-\zeta^{k}\big|^{2}\Big].
        \end{align}
    \end{small}Thus
    \begin{small}
    \[
        \lim_{k\to\infty}\widehat{\mathbb{E}}\Biggl[\Biggl|Y_{t}^{t,\zeta,u}-\sum_{i=1}^{N_{k}}I_{A_{i}^{k}}Y_{t}^{t,x_{i}^{k},u}\Biggr|^{2}\Biggr]=0.
    \]
    \end{small}Moreover, we know that
    \begin{small}
    \[
        V(t,\zeta^{k})=\sum_{i=1}^{N_{k}}I_{A_{i}^{k}}V(t,x_{i}^{k})\leq\sum
        _{i=1}^{N_{k}}I_{A_{i}^{k}}Y_{t}^{t,x_{i}^{k},u}\; \text{q.s.}.
    \]
    \end{small}Hence, we get 
    \begin{small}
    \[
        V(t,\zeta)\leq Y_{t}^{t,\zeta,u}\; \text{q.s.}.
    \]
    \end{small}

    In the following, we prove that for a given 
    $\eta\in L_{G}^{2}(\Omega_{t})$, if $\forall u\in\mathcal{U}[t,T]$, 
    $\eta\leq Y_{t}^{t,\zeta,u}\;$q.s., then $\eta\leq V(t,\zeta)\;$q.s..

    By \eqref{Y-zeta-Y-zetak}, for any $u\in\mathcal{U}[t,T]$, we can obtian
    \begin{small}
    \[
        \eta\leq\sum_{i=1}^{N_{k}}I_{A_{i}^{k}}Y_{t}^{t,x_{i}^{k},u}+\sqrt
        {C\widehat{\mathbb{E}}_{t}\Big[\big|\zeta-\zeta^{k}\big|^{2}\Big]}\; \; \text{q.s.}.
    \]
    \end{small}For fixed $N_{k}$, we get that
    \begin{small}
        \begin{align*}
            \eta \leq\sum_{i=1}^{N_{k}}I_{A_{i}^{k}}V(t,x_{i}^{k})+\sqrt{C
            \widehat{\mathbb{E}}_{t}[\mid\xi-\xi^{k}\mid^{2}]}
             =V(t,\xi^{k})+\sqrt
            {C\widehat{\mathbb{E}}_{t}\Big[\big|\zeta-\zeta^{k}\big|^{2}\Big]}\; \;
            \text{q.s.}.
        \end{align*}
    \end{small}By \eqref{zeta-zetak} and \eqref{Vzeta-Vzetak}, there exists a subsequence
    $(\zeta^{k_i})$ of $(\zeta^k)$ such that 
    \begin{small}
        \begin{equation*}
            V(t,\zeta^{k_i})\to V(t,\zeta),\quad 
            \widehat{\mathbb{E}}_{t}\Big[\big|\zeta-\zeta^{k_i}\big|^{2}\Big]\to 0,\quad \text{q.s.,  as }k_i\to\infty.
        \end{equation*}
    \end{small}Therefore, 
    we deduce that $\eta\leq V(t,\zeta)$ q.s.. 
    The proof is complete.
\end{proof}

Now, let us move on to the derivation of dynamic programming principle.
In order to do this, we introduce the so-called backward semigroup
that was first putted forward by Peng \cite{P97}.

For each initial data $(t,x)$, a positive real number $\delta\leq T-t$ and
$\eta\in\cup_{\varepsilon>0}L_{G}^{2+\varepsilon}(\Omega_{t+\delta})$, we
define
\begin{small}
    \[
    \mathbb{G}_{t,t+\delta}^{t,x,u}[\eta]:=\tilde{Y}_{t}^{t,x,u},
    \]
\end{small}where 
$(X_{s}^{t,x,u},\tilde{Y}_{s}^{t,x,u},\tilde{Z}_{s}^{t,x,u},\tilde{K}_{s}^{t,x,u})_{t\leq s\leq t+\delta}$ 
is the solution of the following $G$-SDE and $G$-BSDE:
\begin{small}
    \begin{align}\label{fbsde-dpp}
        \left\{
        \begin{array}{llll}
            X_{s}^{t,x,u}=x+\int_{t}^{s}b(r,X_{r}^{t,x,u},u_{r})dr+\int_{t}^{s}h^{ij}(r,X_{r}^{t,x
            ,u},u_{r})d\langle B^{i},B^{j}\rangle_{r}
            +\int_{t}^{s}\sigma(r,X_{r}^{t,x,u}
            ,u_{r})dB_{r},\\
            \tilde{Y}_{s}^{t,x,u}=\eta+ \int_{s}^{t+\delta}f(r,X_{r}^{t,x,u},\tilde{Y}_{r}^{t,x,u},\tilde
            {Z}_{r}^{t,x,u},u_{r})dr
            +\int_{s}^{t+\delta}g^{ij}(r,X_{r}^{t,x,u},\tilde{Y}_{r}^{t,x,u},\tilde
            {Z}_{r}^{t,x,u},u_{r})d\langle B^{i},B^{j}\rangle_{r}\\
            \qquad\qquad-\int_{s}^{t+\delta}\tilde{Z}_{r}^{t,x,u}dB_{r}-(\tilde{K}_{t+\delta}^{t,x,u}-\tilde{K}_{s}^{t,x,u}),\ s\in[t,t+\delta].
        \end{array}\right.
    \end{align}
\end{small}Then we have the following dynamic programming principle.
\begin{theorem}\label{thm-DPP}
    Let Assumptions {\bf (A1)}-{\bf (A2)} and {\bf (B1)}-{\bf (B4)} hold. 
    Then for any $t\leq s\leq T$, 
    $x\in\mathbb{R}^{n}$, we have
    \begin{small}
        \begin{align*}
            V(t,x)=  \underset{u(\cdot)\in\mathcal{U}[t,s]}{\text{ess}\inf}
            \mathbb{G}_{t,s}^{t,x,u}\big[V(s,X_{s}^{t,x,u})\big]
            =\underset{u(\cdot)\in\mathcal{U}^{t}[t,s]}{\inf}\mathbb{G}_{t,s}
            ^{t,x,u}[V(s,X_{s}^{t,x,u})].\label{DPP}
        \end{align*}
    \end{small}    
\end{theorem}
\begin{proof}
    On the one hand, by Theorem \ref{V-zeta}, 
    we have for any $u(\cdot)\in\mathcal{U}[t,T]$,
    \begin{small}
    \[
        Y_{s}^{s,X_{s}^{t,x,u},u}\geq V(s,X_{s}^{t,x,u})\; \text{q.s.,}
    \]
    \end{small}where $Y_{s}^{s,X_{s}^{t,x,u},u}=Y_{s}^{t,x,u}$ is the solution of $G$-BSDE
    (\ref{fbsde-dpp}) at time $s$. Thus, 
    by the comparison Theorem \ref{thm3.9},
    we get
    \begin{small}
    \[
        Y_{t}^{t,x,u}\geq\mathbb{G}_{t,s}^{t,x,u}\big[V(s,X_{s}^{t,x,u})\big]\; \text{q.s.},
    \]
    \end{small}which means
    \begin{small}
    \[
        V(t,x)\geq\underset{u(\cdot)\in\mathcal{U}[t,s]}{\text{ess}\inf}
        \mathbb{G}_{t,s}^{t,x,u}\big[V(s,X_{s}^{t,x,u})\big].
    \]
    \end{small}

    On the other hand,   
    since we have the similar estimates under both Lipschitz 
    and monotonicity conditions for the solutions of $G$-BSDEs, 
    just different coefficients, see Theorem 7 in \cite{HJ17}
    and our Theorem \ref{estimate-state-2}.
    Then, by employing the so-called “implied
    partition” approach, 
    we can arrive at result that is identical to that 
    of Lemma 22 in \cite{HJ17}. Furthermore, 
    we can also obtain the same 
    result as the Lemma 24 in \cite{HJ17}.
     
    Hence, by Theorem \ref{deterministic-V}, we get 
    \begin{small}
    \[
        V(t,x)\leq\underset{u(\cdot)\in\mathbb{U}^{t}[t,s]}{\inf}\mathbb{G}
        _{t,s}^{t,x,u}\big[V(s,X_{s}^{t,x,u})\big]=\underset{u(\cdot)\in\mathcal{U}[t,s]}{\text{ess}\inf}\mathbb{G}_{t,s}
        ^{t,x,u}\big[V(s,X_{s}^{t,x,u})\big].
    \]
    \end{small}Finally, we get the desired result.
\end{proof}
\section{The viscosity solution of HJB equation}\label{sec5}
In this section, we aim to establish the connection 
between the value function 
\eqref{valuefunction0} of our concerned stochastic
recursive control problem and the 
viscosity solution of its corresponding HJB equation,
where the HJB equation has a form
\begin{small}
    \begin{align}
        &  \partial_{t}V(t,x)+H(t,x,V,\partial_{x}
        V,\partial_{xx}^{2}V)=0,\label{hjb}\\
        &  V(T,x)=\Phi(x),\quad\ \ x\in\mathbb{R}^{n},\nonumber
    \end{align}
\end{small}associated with the Hamiltonian on 
$[0,T]\times\mathbb{R}^{n}\times\mathbb{R}\times\mathbb{R}^{n}\times\mathbb{S}_{n}$:
\begin{small}
    \begin{align}\label{Hamilton}
        H(t,x,v,p,A)= & \underset{u\in U}{\inf}\bigl\{G(F(t,x,v,p,A,u))+\langle p,b(t,x,u)\rangle
         +f(t,x,v,\sigma
        ^{T}(t,x,u)p,u)\bigr\},\\
        F^{ij}(t,x,v,p,A,u)= & (\sigma^{T}(t,x,u)A\sigma(t,x,u))_{ij}+2\langle
        p,h^{ij}(t,x,u)\rangle+2g^{ij}(t,x,v,\sigma^{T}(t,x,u)p,u),\notag
    \end{align}
\end{small}$(t,x,v,p,A,u)\in[0,T]\times\mathbb{R}^{n}\times\mathbb{R}\times\mathbb{R}^{n}\times\mathbb{S}_{n}\times U$, 
$G$ is defined by equation (\ref{eq-G}).

First, we recall the definition of the viscosity solution, 
refer to \cite{CIL92} for more details.
For every $v \in C([0, T ]\times \mathbb{R}^n)$, denote by $\mathbb{P}^{2,+}v(t, x)$
the “parabolic superjet” of $v$ at $(t, x)$, 
which refers to the set of triples 
$(a, p, X) \in \mathbb{R}\times \mathbb{R}^n\times\mathbb{S}_{n}$
such that
\begin{small}
    \begin{align*}
        v(s,y)\le v(t,x)+a(s-t)+\langle p,y-x\rangle
        +\frac{1}{2}\langle X(y-x),y-x\rangle+o\Big(|s-t|+|y-x|^2\Big).
    \end{align*}
\end{small}Similarly the “parabolic subjet” of $v$ at $(t, x)$ can be 
defined by $\mathbb{P}^{2,-}v(t, x)=-\mathbb{P}^{2,+}(-v)(t, x)$.
\begin{definition}
    For $v \in C([0, T ]\times \mathbb{R}^n)$,
    \begin{enumerate}
        \item [(i)] $v$ is a viscosity subsolution of \eqref{hjb} on
        $[0, T ]\times \mathbb{R}^n$, if $v(T , x) \le \Phi(x)$ and 
        for all $(t, x) \in [0, T ]\times \mathbb{R}^n$,
        \begin{small}
            \begin{align*}
                 a+\inf_{u\in U}\bigl\{G(F(t,x,v,p,A,u))+\langle p,b(t,x,u)\rangle
                +f(t,x,v,\sigma^{T}(t,x,u)p,u)\bigr\}\ge 0, \text{ for } 
                (a, p, X) \in\mathbb{P}^{2,+}v(t, x).
             \end{align*}
        \end{small}
        \item [(ii)] $v$ is a viscosity supersolution of \eqref{hjb} on
        $[0, T ]\times \mathbb{R}^n$, if $v(T , x) \ge \Phi(x)$ and 
        for all $(t, x) \in [0, T ]\times \mathbb{R}^n$,
        \begin{small}
            \begin{align*}
                a+\inf_{u\in U}\bigl\{G(F(t,x,v,p,A,u))+\langle p,b(t,x,u)\rangle
                +f(t,x,v,\sigma^{T}(t,x,u)p,u)\bigr\}\le 0, \text{ for } 
                (a, p, X) \in\mathbb{P}^{2,-}v(t, x).
             \end{align*}
        \end{small}
        \item[(iii)] $v$ is a viscosity solution of \eqref{hjb}, if it is both a viscosity subsolution and
        viscosity supersolution.
    \end{enumerate}
\end{definition}
Then, we indicate the continuity of the value function.
By Lemma \ref{Vx-Vy}, we note that $V(t,x)$ is 
Lipschitz continuous with respect to $x$.
The following lemma gives the continuity of $V(t,x)$ with respect to $t$.
\begin{lemma}\label{lemma-V-holder-t}
    The value function $V(t,x)$ is $\frac{1}{2}$ H\"{o}lder
    continuous in $t$.
\end{lemma}
\begin{proof}
    There is no essential 
    difference between the classical Lipschitz 
    generator and the non-Lipschitz
    generator in our setting about the estimates of the solutions. 
    Hence, the proof is similar to the one of Lemma 25 in \cite{HJ17},
    we leave out here.
\end{proof}

Now, we present the main theorem, that is the existence theorem 
of viscosity solution of HJB equation \eqref{hjb} in this section.
\begin{theorem}\label{viscosity} 
    Let Assumptions {\bf (A1)}-{\bf (A2)} and {\bf (B1)}-{\bf (B4)} hold. 
    $V$ is the value function defined by (\ref{valuefunction0}). 
    Then $V$ is a unique viscosity
    solution of the second-order partial differential equation
    \eqref{hjb}.
\end{theorem}
To prove Theorem \ref{viscosity}, we make some useful preliminaries.

For $l\ge0$, denote
\begin{small}
    \[
    \Pi_l(y):=
    \begin{cases}
        \frac{|y|\land l}{|y|}y,&\quad y\neq0,\\
        0,&\quad y=0.
    \end{cases}
    \]
\end{small}And for $\phi=f,g^{ij}$, set
\begin{small}
    \[
    \phi^l(t,x,y,z,u):=\phi(t,x,\Pi_l(y),z,u).
    \]   
\end{small}Obviously, $\phi^l$ is linear growth in $x, y, z$, i.e.
\begin{small}
    \[
        |\phi^l(t,x,y,z,u)|\le C(1+|x|+|l|+|z|).
    \] 
\end{small}Then we define a sequence of $\phi_n^l,n\in\mathbb{N}$,
\begin{small}
    \begin{equation}\label{eq-fn-txyzv}
        \phi_n^l(t,x,y,z,v)=\inf_{q\in\mathbb{Q}}\{\phi^l(t,x,q,z,v)+n|y-q|\}.
    \end{equation}
\end{small}Note that $\phi_n^l$ has the following properties, which are same as Lemma \ref{lemma4.1}.
\begin{lemma}\label{lemma-fnu}
    Let $C$ be the constant in {\bf (B4)}. We have
    \begin{enumerate}
        
        \item [(i)] $\phi_n^l$ is well defined for $n\ge C$ and we have
        \begin{small}
            \[
                |\phi_n^l(t,x,y,z,u)|\le C(1+|x|+|l|+|z|).
            \]
        \end{small}    
        \item [(ii)] For any $t\in[0,T]$, $x_1,x_2\in\mathbb{R}^n$, $y\in\mathbb{R}$, $z_1,z_2\in\mathbb{R}^d$, $u\in U[t,T]$, 
        \begin{small}
            \[
                |\phi_n^l(t,x_1,y,z_1,u)-\phi_n^l(t,x_2,y,z_2,u)|\le C(|x_1-x_2|+|z_1-z_2|).
            \]
        \end{small}
        \item [(iii)] For any $t\in [0, T],x\in\mathbb{R}^n,y_1, y_2\in\mathbb{R},z\in\mathbb{R}^d,u\in U[t,T]$,
        \begin{small}
            \[
                |\phi_n^l(t,x,y_1,z,u)-\phi_n^l(t,x,y_2,z,u)|\le n|y_1-y_2|.
            \]
        \end{small}
        \item [(iv)] For any $t\in [0, T],x\in\mathbb{R}^n,y\in\mathbb{R},z\in\mathbb{R}^d,u\in U[t,T]$, the 
        sequence $(\phi_n^l(t,x,y,z,u))_n$ is increasing.
        \item[(v)] If $y_n\to y$ as $n\to\infty$, then $f_n^l(t,x,y_n,z,u)\to f^l(t,x,y,z,u)$.
    \end{enumerate}
\end{lemma}

As a result, we construct the following 
sequences of $G$-BSDEs:
\begin{small}
    \begin{align}
        Y_{s}^{t,x,n,l;u}  &  =\Phi(X_{T}^{t,x,u})+\int_{s}^{T}f_n^l(r,X_{r}^{t,x,u},Y_{r}^{t,x,n,l;u},Z_{r}^{t,x
        ,n,l;u},u_{r})dr\notag\\
        &\quad +\int_{s}^{T}g^{ij,l}_n(r,X_{r}^{t,x,u},Y_{r}^{t,x,n,l;u},Z_{r}^{t,x
        ,n,l;u},u_{r})d\langle B^i,B^j\rangle_r\notag\\
        &\quad  -\int_{s}^{T}Z_{r}^{t,x,n,l;u}dB_{r}-(K_{T}^{t,x,n,l;u}-K_{s}^{t,x,n,l;u}),\ \ s\in[
            t,T].\label{eq-fbsde-n-u}
    \end{align} 
\end{small}and 
\begin{small}
    \begin{align}
        Y_{s}^{t,x,l;u}  &  =\Phi(X_{T}^{t,x,u})+\int_{s}^{T}f^l(r,X_{r}^{t,x,u},Y_{r}^{t,x,l;u},Z_{r}^{t,x
        ,l;u},u_{r})dr\notag\\
        &\quad +\int_{s}^{T}g^{ij,l}(r,X_{r}^{t,x,u},Y_{r}^{t,x,l;u},Z_{r}^{t,x
        ,l;u},u_{r})d\langle B^i,B^j\rangle_r\notag\\
        &\quad  -\int_{s}^{T}Z_{r}^{t,x,l;u}dB_{r}-(K_{T}^{t,x,l;u}-K_{s}^{t,x,l;u}),\ \ s\in[
            t,T].\label{eq-fbsde-l-u}
    \end{align} 
\end{small}
\begin{remark}
    Obviously, $\phi^l(\cdot,X_{\cdot}^{t,x,u},y,z,u)\in M_G^\beta(0,T)$,
    this means that $G$-BSDE \eqref{eq-fbsde-l-u}
    admits a unique solution by Theorem \ref{our-well-posedness}.
    Furthermore, by Lemma \ref{lemma-fnu} and assumption {\bf (B1)}, 
    the fact $\phi_n^l(\cdot,X_{\cdot}^{t,x,u},y,z,u)\in M_G^\beta(0,T)$ 
    with $\phi=f,g^{ij}$ is a corollary of Theorem 4.16 in \cite{HWZ16}.
    Then by Theorem \ref{thm2.7}, $G$-BSDE \eqref{eq-fbsde-n-u}
    admits a unique solution. 
\end{remark}
With the solutions of equations \eqref{eq-fbsde-n-u} and \eqref{eq-fbsde-l-u} in hand, we can define the corresponding
sequences of value functions:
\begin{small}
    \begin{equation}\label{valuefunction-Vn}
        V_n^l(t,x):=\underset{u\in\mathcal{U}[t,T]}{ess\inf}\ Y_{t}^{t,x,n,l;u}.
    \end{equation} 
\end{small}and 
\begin{small}
    \begin{equation*}\label{valuefunction-Vl}
        V^l(t,x):=\underset{u\in\mathcal{U}[t,T]}{ess\inf}\ Y_{t}^{t,x,l;u}.
    \end{equation*} 
\end{small}Similarly, for each $n\in\mathbb{N}$, the corresponding Hamiltonians has form:
\begin{small}
    \begin{align}\label{Hamilton-n}
        H_n^l(t,x,v,p,A)=&\inf_{u\in U}\bigl\{G(F_n^l(t,x,v,p,A,u))+\langle p,b(t,x,u)\rangle
        +f_n^l(t,x,v,\sigma
        ^{T}(t,x,u)p,u)\bigr\},
    \end{align} 
\end{small}where
\begin{small}
    \[  
        F_n^{ij,l}(t,x,v,p,A,u)=  (\sigma^{T}(t,x,u)A\sigma(t,x,u))_{ij}+2\langle
        p,h^{ij}(t,x,u)\rangle
        +2g^{ij,l}_n(t,x,v,\sigma^{T}(t,x,u)p,u),
    \] 
\end{small}and 
\begin{small}
    \begin{align*}\label{Hamilton-l}
        H^l(t,x,v,p,A)=\inf_{u\in U}\bigl\{G(F^l(t,x,v,p,A,u))+\langle p,b(t,x,u)\rangle 
        +f^l(t,x,v,\sigma
        ^{T}(t,x,u)p,u)\bigr\},
    \end{align*} 
\end{small}where
\begin{small}
    \[
        F^{ij,l}(t,x,v,p,A,u)=(\sigma^{T}(t,x,u)A\sigma(t,x,u))_{ij}+2\langle
            p,h^{ij}(t,x,u)\rangle
            +2g^{ij,l}(t,x,v,\sigma^{T}(t,x,u)p,u).
    \] 
\end{small}Then, for the above settings, we have the following uniform convergence results.
\begin{lemma}\label{fnu-converge-fu}
    Assume {\bf (B1)}-{\bf (B4)} hold. Then the sequence $\phi_n^l$ 
    defined in \eqref{eq-fn-txyzv} converge to $\phi^l$ uniformly
    in every compact subset of
    $[0,T]\times\mathbb{R}^n\times\mathbb{R}\times\mathbb{R}^n\times U$,
    with $\phi=f,g^{ij}$.
\end{lemma}
\begin{proof}
    By Lasry and Lyons in \cite{LL1986} (see the Theorem on page 4 and Remark (iv) on
    page 5),
    we conclude that $\phi_n^l$ converges to $\phi^l$ uniformly in 
    $[0,T]\times\mathbb{R}^n\times\mathbb{R}\times\mathbb{R}^n\times U$,
    also in each compact subset of
    $[0,T]\times\mathbb{R}^n\times\mathbb{R}\times\mathbb{R}^n\times U$.
\end{proof}
\begin{lemma}\label{Ynu-converge-Yu}
    Suppose {\bf (A1)}-{\bf (A2)} and {\bf (B1)}-{\bf (B4)} hold.
    Then for all $u\in U[t,T]$, we have
    \begin{small}
        \begin{equation*}
            \lim_{l\to\infty}\sup_{(t,x)\in [0,T]\times K_1}\widehat{\mathbb{E}}\biggl[
                \biggl|Y_{t}^{t,x,l;u}- Y_{t}^{t,x,u}\biggr|
            \biggr]=0,
        \end{equation*}
    \end{small}where $K_1$ in an arbitrary compact subset of
    $\mathbb{R}^n$,
    $Y_{\cdot}^{t,x,u}$ is the solution of $G$-BSDE \eqref{state-2} at $\zeta=x$,
    and $Y_{\cdot}^{t,x,l;u}$ is 
    the solution of $G$-BSDEs \eqref{eq-fbsde-l-u}.
\end{lemma}
\begin{proof}
    For simplicity, we assume that $g^{ij}=0$.
    With the compactness of $K_1$, it
    can be covered by a finite $\delta$-net $\{x_1,\dots, x_N \}$. 
    Therefore, for any $x \in K_1$, there exists $x_i$ ,
    such that $|x-x_i| < \delta$. Consequently, we have
    \begin{small}
        \begin{align*}
            &\sup_{0\le t\le T}\widehat{\mathbb{E}}\biggl[|Y_{t}^{t,x,l;u}-Y_{t}^{t,x,u}|\biggr]\\
            &\leq \sup_{0\le t\le T}\biggl( \widehat{\mathbb{E}}\biggl[|Y_{t}^{t,x,l;u}-Y_{t}^{t,x_i,l;u}|\biggr]+
            \widehat{\mathbb{E}}\biggl[|Y_{t}^{t,x_i,u}-Y_{t}^{t,x,u}|\biggr]\biggr)
            +\sup_{0\le t\le T}\widehat{\mathbb{E}}\biggl[|Y_{t}^{t,x_i,l;u}-Y_{t}^{t,x_i,u}|\biggr]\\
            &\leq \sup_{0\le t\le T}\biggl( \sup_{|x-y|\le \delta}\widehat{\mathbb{E}}\biggl[|Y_{t}^{t,x,l;u}-Y_{t}^{t,y,l;u}|\biggr]+
            \sup_{|x-y|\le \delta}\widehat{\mathbb{E}}\biggl[|Y_{t}^{t,x,u}-Y_{t}^{t,y,u}|\biggr]\biggr)
            +\sup_{0\le t\le T}\widehat{\mathbb{E}}\biggl[|Y_{t}^{t,x_i,l;u}-Y_{t}^{t,x_i,u}|\biggr]\\
            &\leq \sup_{|x-y|\le \delta}\biggl( \sup_{0\le t\le T}\widehat{\mathbb{E}}\biggl[|Y_{t}^{t,x,l;u}-Y_{t}^{t,y,l;u}|\biggr]+
            \sup_{0\le t\le T}\widehat{\mathbb{E}}\biggl[|Y_{t}^{t,x,u}-Y_{t}^{t,y,u}|\biggr]\biggr)\\
            &+N\cdot \sup_{i\in(1,\dots,N)}\sup_{0\le t\le T}\widehat{\mathbb{E}}\biggl[|Y_{t}^{t,x_i,l;u}-Y_{t}^{t,x_i,u}|\biggr]
        \end{align*}
    \end{small}By Theorem \ref{estimate-state-2}, 
    we get that for any $\epsilon> 0$, there exists $\delta> 0$, such
    that,
    \begin{small}
        \begin{equation}\label{Ynx-Yny+Yx-Yy}
            \sup_{|x-y|\le \delta}\biggl( \sup_{0\le t\le T}\widehat{\mathbb{E}}\biggl[|Y_{t}^{t,x,l;u}-Y_{t}^{t,y,l;u}|\biggr]+
            \sup_{0\le t\le T}\widehat{\mathbb{E}}\biggl[|Y_{t}^{t,x,u}-Y_{t}^{t,y,u}|\biggr]\biggr)<\epsilon.
        \end{equation}
    \end{small}Applying It\^{o}'s formula to $|Y_{t}^{t,x,l;u}-Y_{t}^{t,x,u}|^2$, by BDG inequality, we have
    \begin{small}
        \begin{align*}
            &|Y_{t}^{t,x,l;u}-Y_{t}^{t,x,u}|^2+\int_{t}^{T}|Z_{s}^{t,x,l;u}-Z_{s}^{t,x,u}|^2d\langle B\rangle_s+2(M_T-M_t)+2(J_T-J_t)\\
            &\le2\int_{t}^{T}\big(Y_{s}^{t,x,l;u}-Y_{s}^{t,x,u}\big)
            \Big(f(s,X_s^{t,x,u},\Pi_l(Y_{s}^{t,x,l;u}),Z_{s}^{t,x,l;u},u_s)-
            f(s,X_s^{t,x,u},Y_{s}^{t,x;u},Z_{s}^{t,x;u},u_s)\Big)ds\\
            &\le \Big(\frac{C^2}{\underline{\sigma}^2}+2\mu\Big)\int_{t}^{T}
            |Y_{s}^{t,x,l;u}-Y_{s}^{t,x,u}|^2ds+\int_{t}^{T}|Z_{s}^{t,x,l;u}-Z_{s}^{t,x,u}|^2d\langle B\rangle_s
            +2\int_{t}^{T}\hat{F}^l(s)ds,
        \end{align*}
    \end{small}where $\small M_t=\int_{0}^{t}(Y_{s}^{t,x,l;u}-Y_{s}^{t,x,u})(Z_{s}^{t,x,l;u}-Z_{s}^{t,x,u})dB_s$ and 
    $\small J_t=\int_{0}^{t}\big[(Y_{s}^{t,x,l;u}-Y_{s}^{t,x,u})^+dK_s^{t,x,l;u}+(Y_{s}^{t,x,l;u}-Y_{s}^{t,x,u})^-dK_s^{t,x;u}\big]ds$
    are $G$-martingales, and 
    \begin{small}
        \begin{align*}
            \hat{F}^l(s)=\big(Y_{s}^{t,x,l;u}-\Pi_l(Y_{s}^{t,x,l;u})\big)
            \big(f(s,X_s^{t,x,u},\Pi_l(Y_{s}^{t,x,l;u}),Z_{s}^{t,x;u},u_s)
            -f(s,X_s^{t,x,u},Y_{s}^{t,x;u},Z_{s}^{t,x;u},u_s)\big).
        \end{align*}
    \end{small}Then by Gronwall inequality, we can obtain
    \begin{small}
        \begin{align}
            \widehat{\mathbb{E}}\Biggl[\sup_{t\le s\le T}|Y_{s}^{t,x,l;u}-Y_{s}^{t,x,u}|^2\Biggr]
            \le C\widehat{\mathbb{E}}\Biggl[\int_{t}^{T}\hat{F}^l(s)ds\Biggr].\label{eq-Yl-Y}
        \end{align}
    \end{small}By Proposition \ref{propA.1}, Proposition \ref{propA.2} 
    and Theorem \ref{estimate-state-1}, in addition to the linear growth property of $f^l$, 
    we have
    \begin{small}
        \begin{align*}
            \widehat{\mathbb{E}}\Biggl[\sup_{t\le s\le T}|Y_{s}^{t,x,l;u}|^{2}\Biggr]
            +\widehat{\mathbb{E}}\Biggl[\int_{t}^{T}|Y_{s}^{t,x,l;u}|^{2}ds\Biggr]+
            \widehat{\mathbb{E}}\Biggl[\int_{t}^{T}|Z_{s}^{t,x,l;u}|^{2}ds\Biggr]\le C(1+|x|^{2}):=C_0,
        \end{align*}
    \end{small}$C_0$ does not depend on $l$ or $t$. Hence, we can calculus and get
    \begin{small}
        \begin{align}
            \widehat{\mathbb{E}}\Biggl[\int_{t}^{T}|\hat{F}^l(s)|^2ds\Biggr]\le C(1+|x|^4):=C_1,\label{eq-F-l}
        \end{align}
    \end{small}uniformly with respect to $l$ and  $C_1$ is independent of $t$.

    Denote
    \begin{small}
        \[
        S_1^l:=\Big\{(s,\omega)\in [t,T]\times \Omega\Big||Y_{s}^{t,x,l;u}|\ge l\Big\},
        \]
    \end{small} 
    \begin{small}
        \[
        S_2^l:=\Big\{\omega\in\Omega\Big| \text{ there exists } t \le
        s \le T , \text{ such that }|Y_{s}^{t,x,l;u}|\ge l\Big\}.
        \]
    \end{small}Then by Markov inequality, we obtain
    \begin{small}
        \begin{align}
            &\widehat{\mathbb{E}}\Biggl[\int_{t}^{T}\mathbf{1}_{S_1^l}ds\Biggr]\le 
            \widehat{\mathbb{E}}\Biggl[\int_{t}^{T}\mathbf{1}_{[0,T]\times S_1^l}ds\Biggr]
            \le \int_{t}^{T}\widehat{\mathbb{E}}\Big[\mathbf{1}_{[0,T]\times S_2^l}\Big]ds
            \le \frac{T}{l^2} \widehat{\mathbb{E}}\Biggl[\sup_{t\le s\le T}|Y_{s}^{t,x,l;u}|^{2}\Biggr]
            \le \frac{C_0T}{l^2}\to 0, \text{ as }l\to \infty.\label{eq-1-S}
        \end{align}
    \end{small}According to \eqref{eq-Yl-Y}-\eqref{eq-1-S}, by H\"{o}lder inequality, we have
    \begin{small}
        \begin{align*}
            &\widehat{\mathbb{E}}\Biggl[\sup_{t\le s\le T}|Y_{s}^{t,x,l;u}-Y_{s}^{t,x,u}|^2\Biggr]
            \le C\widehat{\mathbb{E}}\Biggl[\int_{t}^{T}\hat{F}^l(s)ds\Biggr]\\
            &\le C \Biggl(\widehat{\mathbb{E}}\Biggl[\int_{t}^{T}|\hat{F}^l(s)|^2ds\Biggr]\Biggr)^{\frac{1}{2}}
            \Biggl(\widehat{\mathbb{E}}\int_{t}^{T}\mathbf{1}_{S_1^l}ds\Biggr)^{\frac{1}{2}}
            \le C\frac{\sqrt{C_1C_0T}}{l}.
        \end{align*}
    \end{small}Then for any $t\in[0,T]$, we obtain
    \begin{small}
        \begin{align*}
            \lim_{l\to\infty}\Biggl\{\sup_{0\le t\le T}
            \widehat{\mathbb{E}}\Biggl[
                \sup_{t\le s\le T}|Y_{s}^{t,x,l;u}-Y_{s}^{t,x,u}|^2\Biggr]\Biggr\}=0.
        \end{align*}
    \end{small}Therefore, for any $\epsilon>0$, for fixed $x_i $, 
    there exists $l(i)$, such that
    \begin{small}
        \[
        \sup_{0\le t\le T}\widehat{\mathbb{E}}
        \biggl[|Y_{t}^{t,x_i,l;u}-Y_{t}^{t,x_i,u}|\biggr]<\epsilon,\ l>l(i).
        \]
    \end{small}Set $L =\max_{i=1,...,N}l(i)$, Then for any $l>L$, we have 
    \begin{small}
        \begin{equation}\label{Yxi-Yx}
            \sup_{i\in(1,\dots,N)}\sup_{0\le t\le T}\widehat{\mathbb{E}}\biggl[|Y_{t}^{t,x_i,l;u}-Y_{t}^{t,x_i,u}|\biggr]<\epsilon.
        \end{equation}
    \end{small}Finally, by \eqref{Ynx-Yny+Yx-Yy} and \eqref{Yxi-Yx}, we obtain
    \begin{small}
        \begin{align*}
            \sup_{(t,x)\in [0,T]\times K_1}\widehat{\mathbb{E}}\biggl[|Y_{t}^{t,x,l;u}-Y_{t}^{t,x,u}|\biggr]
            \le \sup_{x\in K_1} \sup_{0\le t\le T}\widehat{\mathbb{E}}\biggl[|Y_{t}^{t,x,l;u}-Y_{t}^{t,x,u}|\biggr]
            < (N+1)\epsilon,
        \end{align*}
    \end{small}for $\delta$ small enough and $l>L$.

    The proof is complete.
\end{proof}
\begin{lemma}\label{Vn-converge-V} 
    Suppose {\bf (A1)}-{\bf (A2)} and {\bf (B1)}-{\bf (B4)} hold.
    Then $V_n^l$ converges to $V^l$ and $V^l$ converges to $V$, uniformly in every compact subset of
    $[0,T]\times\mathbb{R}^n$.
\end{lemma}
\begin{proof}
    Since $Y_t^{t,x,n,l;u} $ converges to $Y_t^{t,x,l;u}$ in $S_G^2(0,T)$,
    we know path-wisely that $Y_t^{t,x,n,l;u} $ converges to $Y_t^{t,x,l;u}$.
    Thus $V_n^l-V^l$ converges to 0 as $n\to\infty$ for any $[0,T]\times\mathbb{R}^n$.
    By Lemma \ref{Vx-Vy} and Lemma \ref{lemma-V-holder-t}, we get 
    $V_n^l-V^l$ is continuous in $(t,x)$. Therefore, by Dini theorem,
    $V_n^l-V^l$ uniformly converges to $0$ in any compact subset of
    $[0,T]\times\mathbb{R}^n$.

    For any $\epsilon>0$, and each $(t,x)\in[0,T]\times\mathbb{R}^n$, by definition of 
    the value function $V(t,x)$,
    we can find $u_1\in U[0,T]$ such that $V(t,x)+\epsilon\ge Y_t^{t,x;u_1}.$
    Then we have 
    \begin{small}
        \[
        V^l(t,x)-V(t,x)\le Y_t^{t,x,l;u_1}-Y_t^{t,x;u_1}+\epsilon.
        \]
    \end{small}Furthermore, there exists $u_2\in U[0,T]$, such that $V^l(t,x)+\epsilon\ge Y_t^{t,x,l;u_2},$
    which means 
    \begin{small}
        \[
        V(t,x)-V^l(t,x)\le Y_t^{t,x;u_2}-Y_t^{t,x,l;u_2}+\epsilon.
        \]
    \end{small}Hence, for each $(t,x)\in[0,T]\times\mathbb{R}^n$, we obtain
    \begin{small}
        \[
        |V^l(t,x)-V(t,x)|\le \widehat{\mathbb{E}}\Big[\Big|Y_t^{t,x,l;u_1}-Y_t^{t,x;u_1}\Big|\Big]
        +\widehat{\mathbb{E}}\Big[\Big|Y_t^{t,x,l;u_2}-Y_t^{t,x;u_2}\Big|\Big]+2\epsilon,
        \]
    \end{small}where $u_1,$ and $u_2$ are given admissible controls.
    Then by H\"{o}lder inequality and Lemma \ref{Ynu-converge-Yu}, 
    for any copmact set $K_2\subset [0,T]\times\mathbb{R}^n$, we have 
    \begin{small}
        \[
        \lim_{n\to\infty}\sup_{(t,x)\in K_2}\biggl\{\widehat{\mathbb{E}}\Big[\Big|Y_t^{t,x,l;u_1}-Y_t^{t,x;u_1}\Big|\Big]
        +\widehat{\mathbb{E}}\Big[\Big|Y_t^{t,x,l;u_2}-Y_t^{t,x;u_2}\Big|\Big]\biggr\}=0.
        \] 
    \end{small}Finally, by the arbitrariness of $\epsilon$, we proof the desired result.
\end{proof}
\begin{lemma}\label{Hn-converge-H}
    Suppose {\bf (A1)}-{\bf (A2)} and {\bf (B1)}-{\bf (B4)} hold.
    Then $H_n^l$ converges to $H^l$ and $H^l$ converges to $H$, uniformly in every compact subset of
    $[0,T]\times\mathbb{R}^{n}\times\mathbb{R}\times\mathbb{R}^{n}\times\mathbb{S}_{n}$.
\end{lemma}
\begin{proof}
    Note that   
    \begin{small}
        \begin{align*}
            H_n^l-H^l
            \leq \inf_{u\in U} \big(G(|g_n^{ij,l}-g^{ij,l}|)+|f_n^l-f^l|\big)
            \le \sup_{u\in U} \big(G(|g_n^{ij,l}-g^{ij,l}|)+|f_n^l-f^l|\big),
        \end{align*}
    \end{small}and 
    \begin{small}
        \begin{align*}
            H^l-H_n^l
            \le \inf_{u\in U} \big(G(|g^{ij,l}-g_n^{ij,l}|)+|f^l-f_n^l|\big)
            \le \sup_{u\in U} \big(G(|g^{ij,l}-g_n^{ij,l}|)+|f^l-f_n^l|\big).
        \end{align*}
    \end{small}Based on the above two inequalities, we get 
    \begin{small}
        \begin{align*}
            |H_n^l-H^l|
            \le \sup_{u\in U} \big(G(|g_n^{ij,l}(t,x,v,z,u)-g^{ij,l}(t,x,v,z,u)|)
            +|f_n^l(t,x,v,z,u)-f^l(t,x,v,z,u)|\big).
        \end{align*}
    \end{small}For any compact $K\subset [0,T]\times\mathbb{R}^{n}\times\mathbb{R}\times\mathbb{R}^{n}\times\mathbb{S}_{n}$,
    we have
    \begin{small}
        \begin{align*}
            &\lim_{n\to\infty}\sup_{(t,x,v,p,A)\in K}|H_n^l(t,x,v,p,A)-H^l(t,x,v,p,A)|\\
            &\leq \lim_{n\to\infty}\sup_{(t,x,v,p,A)\in K}\sup_{u\in U}\big(G(|g_n^{ij,l}(t,x,v,z,u)-g^{ij,l}(t,x,v,z,u)|)
            +|f_n^l(t,x,v,z,u)-f^l(t,x,v,z,u)|\big)\\
            &\leq \lim_{n\to\infty}\sup_{(t,x,y,z,u)\in K}\big(G(|g_n^{ij,l}(t,x,y,z,u)-g^{ij,l}(t,x,y,z,u)|)
            +|f_n^l(t,x,y,z,u)-f^l(t,x,y,z,u)|\big)\\
            &=0
        \end{align*}
    \end{small}The last equality is due to Lemma \ref{fnu-converge-fu} and the continuity of $G$.
    Thus, we prove that $H_n^l$ converges to $H^l$ uniformly in $K$.
    Finally, let $l\to\infty$, we can obtian that 
    \begin{small}
        \[
        \sup_{(t,x,v,p,A)\in K}|H^l(t,x,v,p,A)-H(t,x,v,p,A)|\to 0.
        \]
    \end{small}The proof is complete.
\end{proof}

So far, we show some convergence results which are important to prove Theorem \ref{viscosity}.
Besides, we provide the final preparation, the stability property of viscosity solutions 
(see Lemma 6.2 in \cite{FS06} for details of proof).
\begin{proposition}[Stability]\label{prop-Stability}
    Let $V_n$ be a viscosity subsolution (resp. supersolution) to the following PDE
    \begin{small}
        \begin{align*}
            &  -\partial_{t}V_n(t,x)-H_n(t,x,V_n,\partial_{x}
                V_n,\partial_{xx}^{2}V_n)=0, \ (t,x)\in [0,T]\times \mathbb{R}^n,
        \end{align*} 
    \end{small}where $H_n(t,x,v,p,A):[0,T]\times\mathbb{R}^{n}\times\mathbb{R}\times\mathbb{R}^{n}\times\mathbb{S}_{n}\to\mathbb{R}$
    is continuous and satisfies the ellipticity condition
    \begin{small}
        \[
        H_n(t,x,v,p,A)\le H_n(t,x,v,p,B) \text{ whenever } A\le B.
        \]
    \end{small}Assume that $H_n$ and $V_n$ converge to $H$ and $V$, respectively, uniformly in every compact subset of their own
    domains. Then $V$ is a viscosity subsolution (resp. supersolution) of the limit equation
    \begin{small}
        \begin{align*}
            &  -\partial_{t}V(t,x)-H(t,x,V,\partial_{x}
                V,\partial_{xx}^{2}V)=0, \ (t,x)\in [0,T]\times \mathbb{R}^n.
        \end{align*}
    \end{small}
\end{proposition}
\begin{proof}[{\bf Proof of Theorem \ref{viscosity}}]
The uniqueness of viscosity solution of equation \eqref{hjb} can be proved similarly as in Theorem
5.1 in \cite{ZDP20}. We only prove that $V$ is a viscosity solution of equation \eqref{hjb}.
By Lemma \ref{lemma-fnu}, we note that the conditions of Theorem 26
in \cite{HJ17} are satisfied. Thus, we conclude that 
$V_n^l$ is a viscosity solution of the following second-order partial differential equation:
\begin{small}
    \begin{align*}
        &  -\partial_{t}V_n^l(t,x)-H_n^l(t,x,V_n,\partial_{x}
        V_n,\partial_{xx}^{2}V_n)=0,\\
        &  V_n^l(T,x)=\Phi(x),\quad\ \ x\in\mathbb{R}^{n},
    \end{align*} 
\end{small}where $V_n^l$ and $H_n^l$ are defined by \eqref{valuefunction-Vn} and \eqref{Hamilton-n}, respectively.

Simultaneously, Lemma \ref{Vn-converge-V} and Lemma \ref{Hn-converge-H} 
show the uniform convergence of $V_n^l$ to $V$ and $H_n^l$ to $H$ in each compact subset
of their own domains. Furthermore, it is easy to verify that $H_n^l$ 
fulfills the ellipticity condition. 

Hence, by Proposition \ref{prop-Stability}, we can obtain that $V$
is a viscosity solution of the following limit equation
\begin{small}
    \begin{align*}
        &  -\partial_{t}V(t,x)-H(t,x,V,\partial_{x}
        V,\partial_{xx}^{2}V)=0,\\
        &  V_n(T,x)=\Phi(x),\quad\ \ x\in\mathbb{R}^{n},
    \end{align*} 
\end{small}where $V$ and $H$ are defined by \eqref{valuefunction0} and \eqref{Hamilton}, respectively.
\end{proof}

\section{Applications}\label{sec6}
In this section, we give an example to illustrate the
application of our study to utility.

Specially, we start with a market in which two assets are
available for continuous trading: a riskless asset (bond) 
and a risky asset (stock). 
The bond price process is governed by the ODE
\begin{small}
    \begin{equation*}
        dP^0_t=r_tP^0_tdt,\ P^0_0=1,
    \end{equation*} 
\end{small}and the stock price evolves according to the linear $G$-SDE
\begin{small}
    \begin{equation*}
        dP_t=b_tP_tdt+P_tdB_t,\ P_0=p>0,
    \end{equation*} 
\end{small}where $r:[0,T]\to\mathbb{R}$ is the interest rate of the bond and 
$b:[0, T] \to\mathbb{R}$ is the appreciation rate of the
stock, all of which are continuous functions.

The investor chooses a portfolio strategy 
$\pi:\Omega\times [0,T]\to[-1,1]$ the proportion of the wealth 
to invest into the stock
and a consumption plan $c:\Omega\times[0,T]\to[c_1,c_2]$, $0\le c_1<c_2$.
Then the wealth process $W_t$ of the investor can be described by
\begin{small}
    \begin{equation}\label{Wealth}
        \left\{
        \begin{array}{lll}
            dW_t=\big[r_tW_t+(b_t-r_t)\pi_tW_t-c_t\big]dt+W_t\pi_tdB_t,\\
            W_0=x,
        \end{array}\right.
    \end{equation} 
\end{small}where $x > 0$ is the initial wealth of the investor. 
It is obviously that the state equation \eqref{Wealth} satisfies
assumptions {\bf (A1)} and {\bf (A2)}.

As in classical cases, 
the continuous time Epstein-Zin type 
utility process under volatility uncertainty, we called 
the robust stochastic differential 
utility preference of the investor, which can be characterized 
as the unique solution to the following $G$-BSDE:
\begin{small}
    \begin{equation}\label{utility}
        \left\{
        \begin{array}{lll}
            dU_t=\frac{\delta}{1-\frac{1}{\psi}}(1-\gamma)U_t
            \Biggl[\biggl(\frac{c_t}{((1-\gamma)U_t)^{\frac{1}{1-\gamma}}}\biggr)^{1-\frac{1}{\psi}}-1\Biggr]dt-Z_tdB_t-dK_t,\\
            U_T=h(W_T),
        \end{array}\right.
    \end{equation} 
\end{small}where $h:\mathbb{R}\to\mathbb{R}$ is a given Lipschitz continuous function. 

\begin{remark}
    Volatility uncertainty is particularly important in robust utility and 
    asset pricing models because investors are often 
    uncertain about the true volatility of assets, 
    and this uncertainty can significantly affect their decisions.
    Hence, the process $Z$ is introduced to model volatility uncertainty.
    The fact that $K$ is a decreasing martingale suggests that 
    it can be interpreted as capturing unmodeled risks 
    or ambiguities that are not covered by the $G$-Brownian motion.     
\end{remark}

The optimization objective of the investor is to minimize her utility:
\begin{small}
\[
    \min_{(\pi,c)\in U}U_0,
\] 
\end{small}where
\begin{small}
    \[
        U=\big\{(\pi,c)|(\pi,c):
        [0,T]\times\Omega\to[-1,1]
        \times[c_1,c_2] \text{ and } (\pi,c)\in M_G^2(0,T)\big\} 
    \] 
\end{small}is the admissible control set.
\begin{remark}
    We know that investors typically care about utility maximization, 
    but in our example, utility can be viewed as a form of loss in some senses.
\end{remark}

In general, the aggregator in \eqref{utility} does not 
satisfy the Lipschitz condition with respect to the utility
and the consumption. There are four cases which 
the aggregator is monotonic with respect 
to the utility (see Proposition 3.2 in \cite{KSS13}).
Corresponding to our study, and taking into account the linear growth 
condition {\bf (B4)} with respect to the utility, we focus on two cases:
\begin{small}
    \begin{align*}
        &(i)\gamma>1 \qquad\text{  and   }\qquad \psi>1;\\
        &(ii)\gamma<1 \qquad\text{  and   }\qquad \psi<1.
    \end{align*} 
\end{small}Then we adjust the appropriate power such that 
the generator of \eqref{utility} is continuous 
and monotonic but non-Lipschitz in $\mathbb{R}$ 
with respect to the utility $U$ in the above two cases. 
As to the continuity with respect to the consumption,
both cases are Lipschitz continuous obviously whenever $c_1 > 0$. 
However, if $c_1 = 0$, only case $(i)$ satisfies the
continuous but not Lipschitz continuous condition 
with respect to the consumption.

Therefore, for all suitable non-Lipschitz settings 
which satisfy assumptions {\bf (B1)}-{\bf (B4)}, by 
Theorem \ref{viscosity}, we can conclude that the value function of 
the investor is a viscosity solution of the following HJB equation:
\begin{small}
    \begin{equation*}
        \left\{
        \begin{array}{lll}
            \underset{(\pi,c)\in[-1,1]\times[c_1,c_2]}{\min}\Biggl\{
                -\partial_{t}V(t,x)-\big[x(r_t+\pi(b_t-r_t))-c\big]\partial_{x}V(t,x)
                -G(x^2\pi^2\partial_{xx}^2V(t,x))\\
                \qquad\qquad\qquad-\frac{\delta}{1-\frac{1}{\psi}}(1-\gamma)V(t,x)
                \Biggl[\biggl(\frac{c}{((1-\gamma)V(t,x))^{\frac{1}{1-\gamma}}}\biggr)^{1-\frac{1}{\psi}}-1\Biggr]
            \Biggr\}=0,\\
            V(T,x)=h(x).    
        \end{array}\right.
    \end{equation*} 
\end{small}


\appendix
\renewcommand\thesection{\normalsize Appendix: Some estimates}
\section{ }

\renewcommand\thesection{A}

In this appendix, we give some estimates for the solution of $G$-BSDE:
\begin{small}
    \begin{equation}
        Y_{t} =\xi
        +\int_{t}^{T}f(s,Y_{s},Z_{s})ds
        -\int_{t}^{T}Z_{s}dB_{s}-(K_{T}-K_{t}),\label{eqbsde}.
    \end{equation} 
\end{small}Similar estimates hold for \eqref{GBSDE}.
\begin{proposition}\label{propA.1}
    Let $\xi\in L_{G}^{\beta}(\Omega_T)$ for some $\beta>2$, and 
    $f$ satisfy assumptions 
    {\bf (B1)}, {\bf (B2)}, {\bf ($\text{B3}^\prime$)} and {\bf (B4)}.
    Assume that $(Y, Z, K)\in \mathfrak{S}_{G}^{\alpha}(0,T)$ 
    for some $1<\alpha<\beta$ is a solution of the G-BSDE \eqref{eqbsde}.
    Then there exists a constant 
        $C:=C(\alpha,T,\underline{\sigma},C)>0$ such that
        \begin{small}
            \begin{equation*}
                |Y_{t}|^{\alpha}\leq C\widehat{\mathbb{E}}_{t}\Biggl[|\xi|^{\alpha}+\int
                _{t}^{T}|f_{s}^{0}|^{\alpha}ds\Biggr],
            \end{equation*}
        \end{small}where $f_{s}^{0}=|f(s,0,0)|$.
\end{proposition}       
\begin{proof}
    For any $\gamma,\epsilon>0$, set $\tilde{Y}_{t}=|Y_{t}|^{2}+\epsilon_{\alpha} $, 
    where $\epsilon_{\alpha}=\epsilon(1-\alpha/2)^{+}$, applying It\^{o}'s formula 
    to $\tilde{Y}_{t}^{\alpha/2}e^{\gamma t}$. Then 
    we substitute the monotonic condition {\bf ($\text{B3}^\prime$)} for the global
    Lipschitz condition in the standard deduction  
    to obtain the same forms of estimates as 
    Proposition 3.7 in \cite{HJPS1}.
\end{proof}

\begin{proposition}\label{propA.2}
    Let $\xi\in L_{G}^{\beta}(\Omega_T)$ for some $\beta>2$, and 
    $f$ satisfy assumptions 
    {\bf (B1)}, {\bf (B2)}, {\bf ($\text{B3}^\prime$)} and {\bf (B4)}.
    Assume that $(Y, Z, K)\in \mathfrak{S}_{G}^{\alpha}(0,T)$ 
    for some $1<\alpha<\beta$ is a solution of the G-BSDE \eqref{eqbsde}.
    Then there exists a constant $C=C(\alpha,T,C,\underline{\sigma})>0$
    such that 
    \begin{small}
        \begin{equation}
            \widehat{\mathbb{E}}\Biggl[\biggl(\int_{0}^{T}|Z_{s}|^{2}ds\biggr)^{\frac{\alpha}{2}}\Biggr]
            +\widehat{\mathbb{E}}\Big[|K_{T}|^{\alpha}\Big]\leq C\Biggl\{1+
            \Vert Y\Vert _{\mathbb{S}
            ^{\alpha}}^{\alpha}+
            \widehat{\mathbb{E}}\Biggl[\biggl(\int_{0}^{T}f_{s}^{0}ds\biggr)^{\alpha}\Biggr]\Biggr\}.\label{eq-estimate-zk-alpha}
        \end{equation}
    \end{small} 
\end{proposition}

\begin{proof}
    Applying It\^{o}'s formula to $|Y_{t}|^{2}$, we have
    \begin{small}
        \begin{equation*}
            |Y_{0}|^{2}+\int_{0}^{T}|Z_{s}|^{2}d\langle B\rangle_{s}=|\xi|^{2}+\int
            _{0}^{T}2Y_{s}f(s)ds-\int_{0}^{T}2Y_{s}Z_{s}dB_{s}-\int_{0}^{T}2Y_{s}dK_{s},
        \end{equation*}
    \end{small}where $f(s)=f(s,Y_{s},Z_{s})$. Then
    \begin{small}
        \begin{equation*}
            \Biggl(\int_{0}^{T}|Z_{s}|^{2}d\langle B\rangle_{s}\Biggr)^{\frac{\alpha}{2}}\leq
            C_{\alpha}\Biggl\{|\xi|^{\alpha}+\Biggl|\int_{0}^{T}Y_{s}f(s)ds\Biggr|^{\frac{\alpha}{2}
            }+\Biggl|\int_{0}^{T}Y_{s}Z_{s}dB_{s}\Biggr|^{\frac{\alpha}{2}}+\Biggl|
            \int_{0}^{T}Y_{s}dK_{s}\Biggr|^{\frac{\alpha}{2}}\Biggr\}.
        \end{equation*}
    \end{small}
    
    Replacing the monotonic condition {\bf ($\text{B3}^\prime$)} with the global
    Lipschitz condition, then by BDG inequality and simple calculation, there exists a constant
    $C$ which depends on $C,\alpha,T,\underline{\sigma}$ such that
    \begin{small}
        \begin{align}
            \widehat{\mathbb{E}}\Biggl[\biggl(\int_{0}^{T}|Z_{s}|^{2}ds\biggr)^{\frac{\alpha}{2}}\Biggr]\leq 
            C\Biggl\{ \Vert Y\Vert_{\mathbb{S}^{\alpha}}^{\alpha}+\Vert Y\Vert _{
            \mathbb{S}^{\alpha}}^{\frac{\alpha}{2}}\Biggl[\Big(\widehat{\mathbb{E}}\Big[|K_{T}|^{\alpha
            }\Big]\Big)^{\frac{1}{2}}
            +\Biggl(\widehat{\mathbb{E}}\Biggl[\biggl(\int_{0}^{T}f_{s}^{0}ds\biggr)^{\alpha}\Biggr]\Biggr)^{
            \frac{1}{2}}\Biggr]\Biggr\}. \label{eq-estimate-z-alpha} 
        \end{align}
    \end{small}Moreover,
    \begin{small}
        \begin{equation*}
            K_{T}=\xi-Y_{0}+\int_{0}^{T}f(s)ds-\int_{0}^{T}Z_{s}dB_{s}.
        \end{equation*}
    \end{small}Note that 
    \begin{small}
        \begin{align*}
            |f(s,Y_s,Z_s)|\leq|f(s,Y_s,Z_s)-f(s,Y_s,0)|+|f(s,Y_s,0)|
            \leq C|Z_s|+|f(s,0,0)|+C(1+|Y_s|)
        \end{align*}
    \end{small}By calculation, we obtain
    \begin{small}
        \begin{equation}
            \widehat{\mathbb{E}}\Big[|K_{T}|^{\alpha}\Big]\leq C\Biggl\{1+ \Vert Y\Vert _{\mathbb{S}
            ^{\alpha}}^{\alpha}+\widehat{\mathbb{E}}\Biggl[\biggl(\int_{0}^{T}|Z_{s}|^{2}ds\biggr)^{\alpha/2}\Biggr]+
            \widehat{\mathbb{E}}\Biggl[\biggl(\int_{0}^{T}f_{s}^{0}ds\biggr)^{\alpha}\Biggr]\Biggr\}. \label{eq-estimate-k-alpha} 
        \end{equation}
    \end{small}By \eqref{eq-estimate-z-alpha} and \eqref{eq-estimate-k-alpha}, it is easy to get \eqref{eq-estimate-zk-alpha}.
\end{proof}

\begin{proposition}\label{propA.3}
    Let $\xi^l\in L_{G}^{\beta}(\Omega_T)$ with $\beta>2$, $l=1,2$, and 
    $f^l$ satisfy assumptions 
    {\bf (B1)}, {\bf (B2)}, {\bf ($\text{B3}^\prime$)} and {\bf (B4)}.
    Assume that $(Y^l, Z^l, K^l)\in \mathfrak{S}_{G}^{\alpha}(0,T)$ 
    for some $1<\alpha<\beta$ are solutions of the G-BSDE \eqref{eqbsde}
    corresponding to $\xi^l$ and $f^l$. Set $\hat{Y}_{t}=Y_{t}^{1}-Y_{t}^{2},\hat{Z}
    _{t}=Z_{t}^{1}-Z_{t}^{2}$ and $\hat{K}_{t}=K_{t}^{1}-K_{t}^{2}$.
    Then there exists a constant 
    $C:=C(\alpha,T,\underline{\sigma},C)>0$ such that
    \begin{small}
        \begin{equation}
            |\hat{Y}_{t}|^{\alpha}\leq C\widehat{\mathbb{E}}_{t}\Biggl[|\hat{\xi}
            |^{\alpha }+\int_{t}^{T}|\hat{f}_{s}|^{\alpha}ds\Biggr],  \label{eq-hat-Y}
        \end{equation}
    \end{small}and
    \begin{small}
        \begin{align}
            \widehat{\mathbb{E}}\Biggl[\biggl(\int_{0}^{T}|\hat{Z}_{s}|^{2}ds\biggr)^{\frac{\alpha}{2}}\Biggr]\leq
            C\Biggl\{ \Vert \hat{Y}\Vert_{\mathbb{S}^{\alpha}}^{\alpha}+\Vert \hat {Y}
            \Vert_{\mathbb{S}^{\alpha}}^{\frac{\alpha}{2}}\Biggl\{1+||Y^{1}||_{
            \mathbb{S}^{\alpha}}^{\frac{\alpha}{2}}+||Y^{2}||_{
                \mathbb{S}^{\alpha}}^{\frac{\alpha}{2}}\notag\\
            +\Biggl(\widehat{\mathbb{E}}\Biggl[
            \biggl(\int_{0}^{T}f_{s}^{1,0}ds\biggr)^{\alpha}
            \Biggr]\Biggr)^{\frac{1}{2}}+\Biggl(\widehat{\mathbb{E}}\Biggl[
                \biggl(\int_{0}^{T}f_{s}^{2,0}ds\biggr)^{\alpha}
                \Biggr]\Biggr)^{\frac{1}{2}}\Biggr\}\Biggr\},\label{eq-hat-Z}
        \end{align}
    \end{small}where 
    $\hat{f}_{s}=|f^1(s,Y_{s}^{2},Z_{s}^{2})-f^2(s,Y_{s}^{2},Z_{s}^{2})|$,
    $f_s^{l,0}=|f^l(s,0,0)|$.
\end{proposition}

\begin{proof}
    The proof of \eqref{eq-hat-Y} is similar to the  
    Proposition 3.9 in \cite{HJPS1}, we only need to replace the Lipschitz 
    condition with monotonic condition. Applying It\^{o}'s formula to $|\hat{Y}_{t}|^{2}$, 
    by a similar techniques as that in Proposition \ref{propA.2}, we get 
    \begin{small}
        \begin{align*}
            \widehat{\mathbb{E}}\Biggl[\biggl(\int_{0}^{T}|\hat{Z}_{s}|^{2}ds\biggr)^{\frac{\alpha}{2}}\Biggr]
            \leq C \Biggl\{\Vert \hat{Y}\Vert_{\mathbb{S}^{\alpha}}^{\alpha}+\Vert \hat {Y}
            \Vert_{\mathbb{S}^{\alpha}}^{\frac{\alpha}{2}}\Biggl[
                \biggl(\widehat{\mathbb{E}}\Big[|K_{T}^1|^{\alpha}\Big]\biggr)^{\frac{1}{2}}
                +\biggl(\widehat{\mathbb{E}}\Big[|K_{T}^2|^{\alpha}\Big]\biggr)^{\frac{1}{2}}
                +\biggl(\widehat{\mathbb{E}}\biggl[\int_{0}^{T}|\hat{f}_{s}|^{\alpha}ds\biggr]\biggr)^{\frac{1}{2}}\Biggr]\Biggr\}.
        \end{align*}
    \end{small}By Proposition \ref{propA.2}, we have 
    \begin{small}
        \begin{align*}
            &\biggl(\widehat{\mathbb{E}}\Big[|K_{T}|^{\alpha}\Big]\biggr)^{\frac{1}{2}}
                +\biggl(\widehat{\mathbb{E}}\Big[|K_{T}^\prime|^{\alpha}\Big]\biggr)^{\frac{1}{2}}
                +\biggl(\widehat{\mathbb{E}}\biggl[\int_{0}^{T}|\hat{f}_{s}|^{\alpha}ds\biggr]\biggr)^{\frac{1}{2}}\\
                &\leq C\Biggl\{
                    1+\Vert Y^1
                    \Vert_{\mathbb{S}^{\alpha}}^{\frac{\alpha}{2}}+\Vert Y^2
                    \Vert_{\mathbb{S}^{\alpha}}^{\frac{\alpha}{2}}
                    +\Biggl(\widehat{\mathbb{E}}\Biggl[
                        \biggl(\int_{0}^{T}f_{s}^{1,0}ds\biggr)^{\alpha}
                        \Biggr]\Biggr)^{\frac{1}{2}}+\Biggl(\widehat{\mathbb{E}}\Biggl[
                            \biggl(\int_{0}^{T}f_{s}^{2,0}ds\biggr)^{\alpha}
                            \Biggr]\Biggr)^{\frac{1}{2}}
                \Biggr\}.
        \end{align*}
    \end{small}By the above two inequalities, we get \eqref{eq-hat-Z}.
\end{proof}

\begin{remark}\label{remarkA.4}
    Replacing assumption {\bf ($\text{B3}^\prime$)} with {\bf (B3)}
    in Proposition \ref{propA.1}-\ref{propA.3}, the estimates still hold,
    but $C:=C(\alpha,T,\underline{\sigma},C,\mu)$.
\end{remark}

\section*{Acknowledgement}
We would like to thank Prof. Shige Peng and Dr. Wei He for many helpful discussions.

\section*{Declarations}
We declare that there is no conflict of interest, and 
this manuscript is original, unpublished, and not under review elsewhere.


\begin{thebibliography}{99}

    
    
    
    
        \bibitem{BDHPS03} Briand, Ph., Delyon, B., Hu, Y., Pardoux, E., 
        Stoica, L.: $L^p$ solutions of backward stochastic differential 
        equations. Stoch. Process. Appl. {\bf 108}(1), 109-129 (2003)
    
    
    
        \bibitem{CIL92} Crandall, M.G., Ishii, H., Lions, P.L.: User's guide to 
        viscosity solutions of second order partial differential equations.
        Bull. Am. Math. Soc. {\bf 27}(1), 1-67 (1992)
        
    
        \bibitem {DE} Duffie, D., Epstein, L.: Stochastic differential utility.
        Econometrica {\bf 60}(2), 353-394 (1992)
        
        \bibitem{DHP11} Denis, L., Hu, M., Peng S.: Function spaces
        and capacity related to a sublinear expectation: application to $G$-Brownian
        motion pathes. Potential Anal. {\bf 34}(2), 139-161 (2011)
        
        
    
    
        \bibitem{FS06} Fleming, W.H., Soner, H.M.: Controlled Markov Processes 
        and Viscosity Solutions. Springer Verlag (2006)
    
        \bibitem{He22} He, W.: BSDEs driven by $G$-Brownian motion with non-Lipschitz
        coefficients. J. Math. Anal. Appl. {\bf 505}(2) 125569 (2022)
    
        \bibitem{HJ17} Hu, M., Ji, S.: Dynamic programming principle for stochastic 
        recursive optimal control problem driven by a $G$-Brownian
        motion. Stoch. Process. Appl. {\bf 127}(1), 107-134 (2017)
    
        \bibitem{HJPS1} Hu, M., Ji, S., Peng, S., Song, Y.: Backward stochastic
        differential equations driven by $G$-Brownian motion.
        Stoch. Process. Appl. {\bf 124}(1), 759-784 (2014)
    
        \bibitem{HJb} Hu, M., Ji, S., Peng, S., Song, Y.:Comparison theorem, Feynman-Kac 
        formula and Girsanov transformation for BSDEs driven by $G$-Brownian motion. 
        Stoch. Process. Appl. {\bf 124}(2), 1170-1195 (2014)
    
    
        \bibitem{HQW20} Hu, M., Qu, B., Wang, F.: BSDEs driven 
        by $G$-Brownian motion with time-varying Lipschitz condition. 
        J. Math. Anal. Appl. {\bf 491}(2), 124342 (2020)
    
    
        \bibitem{HW2018} Hu, M., Wang, F.: Stochastic optimal 
        control problem with infinite horizon driven by $G$-Brownian motion.
        ESAIM Control Optim. Calc. Var. {\bf 24}(2), 873-899 (2018)
    
        \bibitem{HJL22} Hu, M., Ji, S., Li, X.: Dynamic programming principle and
        Hamilton-Jacobi-Bellman equation under nonlinear expectation. 
        ESAIM Control Optim. Calc. Var. {\bf 28}(5), 1-21 (2022)
    
        \bibitem{HWZ16} Hu, M., Wang, F., Zheng, G.: Quasi-continuous random variables and processes under the $G$
        expectation framework. Stoch. Process. Appl. {\bf 126}(8), 2367-2387(2016) 
    
        \bibitem{HYPSG95} Hu, Y., Peng, S.: Solution of 
        forward-backward stochastic differential equations. 
        Probab. Theory Relat. Fields {\bf 103}(2), 273-283 (1995)
    
        \bibitem{KS14} Karatzas, I., Shreve, S.: Brownian Motion and Stochastic Calculus. Springer (2014)
    
        \bibitem{KSS13} Kraft, H., Seifried, F.T., Steffensen, M.: Consumption-portfolio optimization 
        with recursive utility in incomplete markets. Finance Stoch. {\bf 17}(1), 161-196 (2013)
    
        
        \bibitem{LL1986} Lasry, J.M., Lions, P.L.: A remark on regularization in Hilbert spaces. 
        Israel J. Math. {\bf 55}(3), 257-266 (1986)
    
        \bibitem{LSM97} Lepeltier, J.P., San Martin, J.: Backward stochastic differential equations
        with continuous coefficient. Stat. Probab. Lett. {\bf 32}(5), 425-430 (1997)
    
        \bibitem{Lipeng09} Li, J., Peng, S.: Stochastic optimization 
        theory of backward stochastic differential equations
        with jumps and viscosity solutions of 
        Hamilton-Jacobi-Bellman equations. Nonlinear Anal. {\bf 70}(4), 1776-1796 (2009)
    
        
        \bibitem{Liu20} Liu, G.: Multi-dimensional BSDEs driven by 
        $G$-Brownian motion and related system of fully nonlinear PDEs. Stochastic 
        {\bf 92}(5), 659-683 (2020)
    
    
        \bibitem{Pardoux99} Pardoux, E.: BSDE's, weak convergence 
        and homogenization of semilinear PDE's, in Nonlinear Analysis, Differential Equations
        and Control. Kluwer Academic Publishers, Dordrecht, 503-549 (1999)
        
    
        \bibitem{PP99} Pardoux, E., Tang, S.: Forward-backward 
        stochastic differential equations and quasilinear 
        parabolic PDEs. Probab. Theory Relat. Fields {\bf 114}(2), 123-150 (1999)
        
        \bibitem{Peng1991} Peng, S.: Probabilistic Interpretation for
        Systems of Quasilinear Parabolic Partial Differential Equations. 
        Stochastics {\bf 37}(1-2), 61-74 (1991)
        
        
        \bibitem{Peng1992} Peng, S.: A Generalized Dynamic Programming
        Principle and Hamilton-Jacobi-Bellmen equation. Stochastics {\bf 38}, 119-134 (1992)
    
        
        \bibitem{P97} Peng, S.: Backward stochastic differential 
        equations-stochastic optimization theory and viscosity solutions of HJB equations.
        In Topics on Stochastic Analysis, 
        edited by Yan, J., Peng, S., Fang, S., Wu, L. Science Press, Beijing (in Chinese),
        85-138 (1997)
    
        
        
        
        \bibitem{P07a} Peng, S.: $G$-expectation, $G$-Brownian Motion and
        Related Stochastic Calculus of It\^o type. Stochastic Analysis and
        Applications, 541-567, Abel Symposium, 2, Springer, Berlin (2007)
        
        
        
        
        
        \bibitem{P10} Peng, S.: Nonlinear Expectations and Stochastic
        Calculus under Uncertainty. arXiv preprint arXiv:1002.4546v1 (2010)
        
    
        \bibitem{P19} Peng, S.: Nonlinear Expectations and Stochastic Calculus
        Under Uncertainty: With Robust CLT and $G$-Brownian Motion. 
        Probability Theory and Stochastic Modelling 95, Springer (2019)
    
    
        \bibitem{POD13} Possama\"{i}, D.: Second order backward stochastic 
        differential equations under a monotonicity condition. 
        Stoch. Process. Appl. {\bf 123}(5), 1521-1545 (2013) 
    
        \bibitem{PZ18} Pu, J., Zhang, Q.: Dynamic programming principle 
        and associated Hamilton-Jacobi-Bellman equation for 
        stochastic recursive control problem with non-Lipschitz aggregator.
        ESAIM Control Optim. Calc. Var. {\bf 24}(1), 355-376 (2018)
    
        
        \bibitem{STZ12} Soner, M., Touzi, N., Zhang, J.: Well-posedness of
        Second Order Backward SDEs. Probab. Theory Relat. Fields
        {\bf 153}(1-2), 149-190 (2012)

        \bibitem{Song2019} Song, Y.: Backward Stochastic Differential Equations Driven by G-Brownian
        Motion Under a Monotonicity Condition. Chinese Annals of Mathematics, Series A. {\bf 40}(2), (2019). 
        (doi:10.16205/j.cnki.cama.2019.0015)
        
        \bibitem{Song11} Song, Y.: Some properties on G-evaluation and
        its applications to G-martingale decomposition. Sci. China Math.
        {\bf 54}(2), 287-300 (2011)
        
        
        \bibitem{WZ21} Wang, F., Zheng, G.: Backward stochastic differential equations 
        driven by $G$-Bsrownian motion with
        uniformly continuous generators. J. Theor. Probab. {\bf 34}(2), 660-681 (2021)
    
        \bibitem{WWZ22} Wen, J., Wu, Z., Zhang, Q.: Dynamic programming 
        principle for delayed stochastic
        recursive optimal control problem and HJB equation
        with non-Lipschitz generator. arXiv Preprint arXiv:2205.03052v3 (2022)
    
        \bibitem{Wuyu08} Wu, Z., Yu, Z.: Dynamic programming principle for 
        one kind of stochastic recursive optimal
        control problem and Hamilton-Jacobi-Bellman equation. 
        SIAM J. Control. Optim. {\bf 47}(5), 2616-2641 (2008)
    
        \bibitem{Zhao24} Zhao, B.: BSDEs driven by G-Brownian motion with time-varying uniformly
        continuous generators. arXiv Preprint arXiv:2409.16574v1 (2024)
        
        \bibitem{ZDP20} Zhuo, Y., Dong, Y., Pu, J.: Dynamic programming 
        principle and viscosity solutions
        of Hamilton-Jacobi-Bellman equations for 
        stochastic recursive control problem 
        with non-Lipschitz generator. 
        Appl. Math. Optim. {\bf 82}(2), 851-887 (2019)
    \end{thebibliography}
\end{document}